\documentclass[12pt,a4paper]{article}
\usepackage{latexsym,amssymb,amsmath,mathrsfs}
\usepackage{sectsty}
\usepackage{color}
\usepackage{bm}
\usepackage[T1]{fontenc}

\definecolor{green1}{rgb}{0,0.392157,0}
\definecolor{blue1}{cmyk}{0.4,0.4,0,0.4}
\usepackage[anchorcolor=black,%
linkcolor=green1,%
colorlinks=true,%
citecolor=blue1,%
bookmarks=true,%
bookmarksnumbered=true,%
]{hyperref}
\usepackage{cite}

\allsectionsfont{\normalsize\sc\center}

\usepackage[top=3cm,bottom=3.5cm,left=2.8cm,right=2.8cm]{geometry}

\newtheorem{Thm}{Theorem}[section]
\newtheorem{Prop}[Thm]{Proposition}

\newtheorem{Lem}[Thm]{Lemma}
\newtheorem{Defn}[Thm]{Definition}
\newtheorem{Rem}[Thm]{Remark}

\makeatletter \@addtoreset{equation}{section} \makeatother

\newcommand{\nn}{\nonumber}
\newcommand{\ds}{\displaystyle}

\newenvironment{Proof}{\par\begin{trivlist}%
\item[]{\bf Proof.\ }}%
{\hfill $\square$ \end{trivlist}\par}
\newenvironment{tProof}[1]{\par\begin{trivlist}%
\item[]{\bf Proof of #1.\ }}%
{\hfill $\square$ \end{trivlist}\par}

\newcommand{\N}{\mathbb{N}}
\newcommand{\R}{\mathbb{R}}
\renewcommand{\d}{\mathrm{d}}
\DeclareMathOperator{\supp}{supp}
\DeclareMathOperator{\Ric}{Ric}
\DeclareMathOperator{\Hess}{Hess}
\DeclareMathOperator{\Ent}{Ent}
\newcommand{\lip}{\mathrm{lip}}
\newcommand{\Lip}{\mathrm{Lip}}
\newcommand{\e}{\mathrm{e}}
\newcommand{\m}{\mathfrak{m}}

\title{\sf Monotonicity and rigidity of the $\mathcal{W}$-entropy
on $\mathsf{RCD} (0,N)$ spaces}
\author{Kazumasa Kuwada\thanks{Mathematical Institute, Tohoku University,
Sendai 980-8578, Japan ({\sf kuwada@m.tohoku.ac.jp}). Supported in part by JSPS Grant-in-Aid for Young Scientist (A) (KAKENHI) 26707004}\hspace{.5em} and Xiang-Dong Li\thanks{Academy of Mathematics and Systems Science, Chinese Academy of Sciences, 55, Zhongguancun East Road, Beijing, 100190, China
({\sf xdli@amt.ac.cn})
and School of Mathematical Sciences, University of Chinese Academy of Sciences, Beijing, 100049, China. Research supported by NSFC No. 11771430, Key Laboratory RCSDS, CAS, No. 2008DP173182, and Hua Luo-Keng Research Grant of AMSS, CAS (2015-2017).
}}

\date{}
\begin{document}
\maketitle
\begin{abstract}
By means of a space-time Wasserstein control,
we show the monotonicity of the $\mathcal{W}$-entropy functional
in time along heat flows on possibly singular metric measure spaces
with non-negative Ricci curvature
and a finite upper bound of dimension in an appropriate sense.
The associated rigidity result on the rate of dissipation of
the $\mathcal{W}$-entropy is also proved.
These extend known results even on weighted Riemannian manifolds
in some respects.
In addition, we reveal that some singular spaces will exhibit the rigidity models
while only the Euclidean space does in the class of
smooth weighted Riemannian manifolds.
\end{abstract}


\section{Introduction}
\label{sec:intro}
The $\mathcal{W}$-entropy functional was first
introduced on Ricci flow in the celebrating work of G.~Perelman \cite{Perel_RF1}
for the resolution of the Poincar\'e conjecture.
Since then, it has played important role in the various topics in the study of
geometric analysis and stochastic analysis. Among others, the $\mathcal{W}$-entropy exhibits two specific properties:
monotonicity in time along conjugate heat equations and rigidity.
The latter means that the time derivative of the $\mathcal{W}$-entropy along
a conjugate heat equation vanishes if and only if
the Ricci flow is a gradient shrinking Ricci soliton.
As one of possible extensions in connection with these properties,
the notion of $\mathcal{W}$-entropy is exported
in several different situations.
L.~Ni \cite{Ni_W-ent0,Ni_W-ent1} brought
the notion of $\mathcal{W}$-entropy to the heat equation $\partial_t u=\Delta u$ on \emph{static} Riemannian manifolds,
where static means that the Riemannian metric does not depend on time.
He proved the same sort of monotonicity and rigidity for the $\mathcal{W}$-entropy for the heat equation $\partial_t u=\Delta u$
under non-negative Ricci curvature with an additional bounded geometry assumption. 
In a series of works by the second named author \cite{LiXD_W-ent1, LiXD_W-ent0, LiXD_W-ent2, LiXD_W-ent3}, the $\mathcal{W}$-entropy formula as well as its monotonicity and rigidity have been extended to the heat equation $\partial_t u=Lu$ associated with the Witten Laplacian $L=\Delta-\nabla\phi\cdot\nabla$ on complete Riemannian manifolds with weighted volume measure $d\mu=\e^{-\phi}dv$, where $v$ is the Riemannian volume measure. Furthermore, such extensions have been carried out for the heat equation associated with the time-dependent Witten Laplacian on the so-called $(K, m)$ super-Ricci flows by the second named author and S.~Li \cite{LiLi_W-ent1, LiLi-JFA2018, LiLi-SCM2018, LiLi-AJM2018}. 
Here the $(K, m)$ super-Ricci flow means that the Riemannian metric $g_t$
is time-dependent and evolves along the following differential inequality
 \[
{1\over 2} {\partial g_t\over \partial t} +\Ric_{m, n}(L)\geq Kg_t,
 \]
 where $\Ric_{m, n}(L)=\Ric(g_t)+\nabla^2\phi_t-{\nabla\phi_t\otimes \nabla\phi_t\over m-n}$ is the $m$-dimensional Bakry-Emery Ricci curvature on $n$-dimensional complete Riemannian manifolds $(M, g_t)$ with fixed weighted volume measure $d\mu=e^{-\phi_t}dv_{g_t}$.
Note that the notion of the $(K, m)$-super Ricci flow can be regarded as
an extension of the super Ricci flow in geometric analysis 
\[
 {\partial g_t\over \partial t}\geq  -2\Ric,
 \]
which includes the Ricci flows as the case of equality, and also an extension of static Riemannian manifolds
with Ricci curvature bounded from below by a constant (i.e., $\Ric\geq Kg$).

In this article, we study the same sort of problem on more singular spaces
than differentiable manifolds. The notion of spaces with a lower Ricci curvature
bound has been extended from Riemannian manifolds to metric measure spaces
by means of optimal transport \cite{Lott-Vill_AnnMath09,Sturm_Ric12}.
Thus it seems natural to consider this problem on such spaces
with non-negative Ricci curvature in a generalized sense.
As an initial work to this direction, we only consider the static case,
though it seems possible to consider the corresponding problem 
on time-dependent metric measure spaces on the basis of \cite{Kopf-St,Sturm_Rf1}.

To state our result in comparison with the one in the smooth case,
let us begin with reviewing the result on $n$-dimensional 
weighted Riemannian manifolds as mentioned above 
according to \cite{LiXD_W-ent0,LiXD_W-ent2, LiLi_W-ent1, LiLi-JFA2018, LiLi-SCM2018}.
To begin with, let us remark that, for $K \in \R$, 
Bakry-\'Emery's curvature-dimension condition for $L$
\begin{equation} \label{eq:intro-BE}
\frac12 L | \nabla f |^2 - \langle \nabla f , \nabla L f \rangle
\ge
K | \nabla f |^2 + \frac{1}{m} | L f |^2
\end{equation}
holds for any $f \in C^3 (M)$ if and only if
\begin{equation} \label{eq:wRic}
\Ric_{m,n} \ge K.
\end{equation}
Roughly speaking, $m$ plays the role of
an upper bound of the dimension of the space.
Indeed, \eqref{eq:intro-BE} can be used as an abstract generalization of
the condition ``$\Ric \ge K$ and $\dim \le m$'' in Bakry-\'Emery theory
(see \cite{BE:diff-hyp,BGL}).
Suppose $n < m < \infty$.
Then we define the $\mathcal{W}$-entropy 
for $f \in C^1 (M)$ and $t > 0$ as follows:
\begin{equation} \label{eq:intro-W}
\mathcal{W} ( f, t )
 : =
\int_M
\left[
    t | \nabla f |^2 + m - f
\right]
\frac{\e^{-f}}{ ( 4 \pi t )^{m/2} }\, \d \mathfrak{m}.
\end{equation}
It coincides with the one introduced in \cite{Ni_W-ent0} when $m = n$ and $\phi = 0$.
Suppose that $f$ depends also on $t$ and that
$u : = \e^{-f} / ( 4 \pi t )^{m/2}$ solves the heat equation $\partial_t u = L u$.
Then, the monotonicity of $\mathcal{W}$
\begin{equation} \label{eq:intro-mono}
\dfrac{\d}{\d t} \mathcal{W} ( f , t ) \le 0
\end{equation}
holds under $\Ric_{m.n} \ge 0$ and a bounded geometry assumption.
Suppose additionally that $u$ is a heat kernel.
That is, $u \to \delta_y$ for some $y\in M$ weakly.
Then the rigidity of $\mathcal{W}$ states that the equality holds
in \eqref{eq:intro-mono} at some $t > 0$
only when $(M,g)$ is isometric to $\R^m$ and $\phi$ is constant.
The proof is based on establishing the so-called ``entropy formula''
which explicitly describes the time derivative of $\mathcal{W}$. 
For details, see \cite{LiXD_W-ent0,LiXD_W-ent2, LiLi_W-ent1, LiLi-JFA2018, LiLi-SCM2018}. 

A natural class of metric measure spaces where we consider our problem is
${\sf RCD}(0,N)$ (or equivalently ${\sf RCD}^* (0,N)$) spaces.
Intuitively, ${\sf RCD}^* (K,N)$ means that
the space satisfies ``$\Ric \ge K$ and $\dim \le N$''
and the canonical heat flow given by the metric measure structure
is linear in initial data. There are several different characterizations
of ${\sf RCD}^* (K,N)$ spaces, and the Bochner inequality like \eqref{eq:intro-BE}
is one of them. See \cite{Ambrosio:2012tp,AGS3,AGS_BE-CD} for $N = \infty$ and
\cite{Erbar:2013wf,AMS} for $N < \infty$.
Note that the study of ${\sf RCD} (K,N)$ spaces for finite $N$ are
initiated in \cite{Gigli_diff-str,Gigli_Split} and connection with
the Bakry-\'Emery condition is established in \cite{Erbar:2013wf,AMS}.
See the next section for more details and additional references.
Recall that any $n$-dimensional weighted Riemannian manifold $( M , g , \mu )$ 
with $\Ric_{N.n} \ge K$ is an ${\sf RCD}^* (K,N)$ space 
by regarding it as a metric measure space by the Riemannian distance
and the weighted measure $\mu$.

Our main theorems are the monotonicity of $\mathcal{W}$
(Theorem~\ref{thm:mono-W})
and the associated rigidity (Theorem~\ref{thm:rigid-W})
on ${\sf RCD} (0,N)$ spaces.
Note that our rigidity theorem improves the previous result even
on weighted Riemannian manifolds in the following three respects.
First, we do not require a differentiability of $\mathcal{W}$ in time
along the heat flow but we consider the right upper derivative instead.
Second, we do not need to assume the initial data to be the Dirac mass.
Actually, it follows as a consequence of the rigidity:
If the right upper derivative vanishes, then the initial data must be Dirac.
Third, we do not require any assumption corresponding to the ``bounded geometry''.
Moreover, we find that not only Euclidean spaces but the Euclidean cones
enjoys the vanishing time derivative of $\mathcal{W}$, where the vertex of the cone
coincides with the point where the initial mass is located.
The Euclidean cones appearing in our rigidity have singularity
at vertex if it is not a Euclidean space, and in general it is even not a manifold.
In this sense, our result is compatible with the previous ones and
we succeed in finding new examples as a result of expanding
the class of spaces we consider.

For the proof of our main results, we relies on an approach from optimal transport
which is not used in previous results \cite{LiXD_W-ent1,LiXD_W-ent0,LiXD_W-ent3, LiLi_W-ent1}.
A naive approach to our problem is to establish the entropy formula on
${\sf RCD} (0,N)$ spaces. However, it does not seem to be straightforward
by the following two reasons.
First, the entropy formula involves the Ricci curvature $\Ric$, the Hessian $\Hess f$ and some second order tensor calculus is required.
Although such objects have been introduced on ${\sf RCD}$ spaces
in recent development \cite{Gigli_NonDG},
the lack of smoothness can be an obstacle.
Second, it seems that some assumption like the bounded geometry is required
to obtain the entropy formula. It is not clear how we formulate
such an assumption with keeping non-trivial examples since the bounded geometry
assumption involves the Riemann curvature tensor and its derivatives.
Because of them, we prove the monotonicity and the rigidity
\emph{without} entropy formulae.
The idea of our proof comes from \cite{Topp_Lopt},
where P.~Topping studies the monotonicity
of the $\mathcal{W}$-entropy on a Ricci flow on a compact manifold
by proving an estimate of a transportation cost between two heat distributions
where the cost function is given by Perelman's $\mathcal{L}$-distance.
In the static case, his estimate reduces to a so-called space-time
$L^2$-Wasserstein control for heat distributions \eqref{eq:control},
ant it is indeed one of
characterizing properties of ${\sf RCD} (0,N)$ space \cite{Erbar:2013wf,K14}.
On the one hand, heat flow can be regarded as a gradient flow of the relative
entropy functional (or the Boltzmann-Shannon entropy with an opposite sign)
on $L^2$-Wasserstein space, and thus there is a strong connection between
heat flows and the relative entropy by means of optimal transport.
On the other hand, we can write the $\mathcal{W}$-entropy
by using the relative entropy and its dissipation along heat flow
(or the Fisher information)
\cite{LiXD_W-ent0,LiXD_W-ent1,LiXD_W-ent2,LiXD_W-ent3}.
The proof of the monotonicity follows from combining these two observations 
with the idea in \cite{Topp_Lopt}.
Note that the monotonicity of $\mathcal{W}$ is already studied
on ${\sf RCD} (0,N)$ spaces by different means.
This problem is considered first in \cite{Jiang:2015}
when the underlying space is compact.
The noncompact case is discussed in \cite{LiHQ_W-ent}
by following an argument in \cite{BauGaro-W} on Riemannian manifolds.
However, it seems that some technical details are not well described
in the latter case.
The proof of the rigidity also relies on the space-time Wasserstein control
in its first step, but we use the condition in a more subtle way.
It implies an identity for the Fisher information, and the final conclusion
is reduced to the recent result on the volume rigidity \cite{Gig_DePhi}
(See Theorem~\ref{thm:G-DP} below)
by using recently developed analytic tools such as the Li-Yau inequality
\cite{Jiang:2014vua}
and the Varadhan type short time asymptotic for the heat kernel
following from \cite{Jiang:2014wo}.
Roughly speaking, the assumption of the rigidity of the $\mathcal{W}$-entropy
implies the equality in the Laplacian comparison theorem
(See Proposition~\ref{prop:id-d} and Remark~\ref{rem:id-d}). From geometric
viewpoint, it is almost equivalent to the equality in Bishop-Gromov inequality
studied in \cite{Gig_DePhi}.
From probabilistic viewpoint, the transition probability becomes
the Gaussian kernel if it starts from the reference point
(see Proposition~\ref{prop:hat-GF}).
A typical, elementary but non-Euclidean example for the latter one
is the $N$-dimensional Bessel process with $N \in [1,\infty)$
starting from the endpoint of the interval.
In this case, we take a weighted measure $\mathfrak{m} (\d r) = r^{N-1} \d r$
to make $[0,\infty)$ equipped with the Euclidean distance
to be a metric measure space.

The structure of this article is as follows.
In the next section, we introduce several notions concerning with
metric measure spaces and ${\sf RCD}$ spaces.
Known properties of ${\sf RCD} (0,N)$ spaces
we will use in the sequel are also prepared there.
The monotonicity of the $\mathcal{W}$-entropy and the rigidity are
shown in Section \ref{sec:mono} and Section \ref{sec:rigid} respectively.
For the rigidity, we first argue that the case when the initial data is Dirac,
after showing that it certainly happens at some point (Lemma~\ref{lem:id-I}).
By using the consequence of it, we show that the initial data must be Dirac
(Lemma~\ref{lem:Dirac}).
Some results related with our main theorem are gathered
in Section~\ref{sec:related}.
We deal with four different topics there.
First, we show that the heat flow becomes also $L^2$-Wasserstein geodesic
if the assumption of the rigidity holds (Proposition~\ref{prop:hat-G}).
Indeed, as already observed by the second named author and S. Li in \cite{LiLi_W-ent3, LiLi-SCM2018}, there is some similarity in the study of $\mathcal{W}$-entropy for the heat flow on the underlying manifolds and the geodesic flow on the $L^2$-Wasserstein space. To understand such a similarity better, the second named author and S. Li  \cite{LiLi_W-ent3, LiLi-SCM2018} introduced the Langevin deformation of flows over weighted  Riemannian manifolds, which can be regarded as a natural interpolation between the heat flow on underlying manifolds and the geodesic flow on the Wasserstein space equipped with Otto's infinite dimensional Riemannian metric. Moreover, Perelman's $\mathcal{W}$-entropy formula has been extended in \cite{LiLi_W-ent3, LiLi-SCM2018} to the geodesic flow and the Langevin deformation of flows on the Wasserstein space on Riemannian manifolds with non-negative Ricci curvature and on weighted Riemannian manifolds with non-negative $m$-dimensional Bakry-Emery Ricci curvature.
From this point of view, it seems meaningful that the same property holds
even in the framework of ${\sf RCD}$ spaces. 
Second, we discuss some relations 
between the (logarithmic) Sobolev inequality 
and the $\mathcal{W}$-entropy in our framework. 
It is already pointed out in \cite{Perel_RF1} that $\mathcal{W}$-entropy 
is related with the logarithmic Sobolev inequality. 
See also \cite{LiXD_W-ent1,LiXD_W-ent0, LiLi_W-ent1, LiLi-JFA2018}. 
Third, we consider a stronger rigidity result under a stronger assumption 
(Theorem~\ref{thm:s-rigid-W}).
In this case, the conclusion becomes the same
as the rigidity theorem on weighted Riemannian manifolds which was proved previously in \cite{LiXD_W-ent1,LiXD_W-ent2, LiLi_W-ent1}.
Fourth, we consider the almost rigidity. Here ``almost rigidity''
asserts that the conclusion of the rigidity almost holds
if the assumption of the rigidity is almost satisfied.
That is, a weaker assumption implies a weaker conclusion.
The most famous almost rigidity result would be
an extension of the Cheeger-Gromoll splitting theorem. See e.g. \cite{FLZ} for such an extension on weighted Riemannian manifolds with non-negative finite dimensional or infinite dimensional Bakry-Emery Ricci curvature.
Now the most general ``almost splitting theorem'' is
formulated in the framework of
${\sf RCD} (0,N)$ spaces \cite{Gigli_Split} (See references therein also).
The key property for the almost rigidity in \cite{Gigli_Split} is
that the ${ \sf RCD } (0,N)$ condition is stable under
a (pointed measured Gromov-Hausdorff) convergence of metric measure spaces.
Though our assumption on the rigidity of the $\mathcal{W}$-entropy seems
less stable under convergence of spaces,
we are somehow able to formulate an almost rigidity.

\section{Framework}
\label{sec:frame}

Let $(X,d)$ be a complete and separable geodesic metric space.
Here ``geodesic'' means that
for any $x_0, x_1 \in X$, there exists $\gamma : [ 0, 1 ] \to X$
such that $\gamma_i = x_i$ ($i=0,1$) and
$d ( \gamma_s , \gamma_t ) = |s - t | d ( x_0 , x_1 )$.
We call such $\gamma$ a (minimal) geodesic joining $x_0$ and $x_1$.
Let $\m$ be $\sigma$-finite Borel measure on $X$.
Suppose that $\m(B_r(x)) \in ( 0 , \infty )$
for any metric ball $B_r (x)$ of radius $r > 0$ centered at $x \in X$.
In particular, $\supp \m = X$ holds.
We call the triplet $(X,d,\m)$ a metric measure space in this article.
A typical example of metric measure space
we should have in mind is the weighted Riemannian manifold
as reviewed in the introduction.

Both for defining $\mathsf{RCD}$ spaces and for considering
the canonical heat flow on $(X,d,\m)$,
we require the notion of the ($L^2$-)Cheeger energy functional.
Let $\Lip (X)$ be the set of all Lipschitz continuous functions on $X$
and $\Lip_b (X) = \Lip (X) \cap L^\infty (\mathfrak{m})$.
For $f \in \Lip (X)$,
we define the local Lipschitz constant $\lip(f)(x)$ of $f$ at $x \in X$ by
\[
\lip(f)(x) : = \varlimsup_{y \to x} \frac{ | f (y) - f (x) | }{ d (x,y) }.
\]
We regard $\lip (f)$ as a function on $X$.
By means of local Lipschitz constant, we define
the Cheeger energy $\mathsf{Ch}$ as follows:
for $f \in L^2 (\m)$,
\[
\mathsf{Ch} (f) :=
\frac12
\inf \left\{ \left.
    \liminf_{n \to \infty}
    \int_X \lip ( f_n )^2 \d \m
    \; \right| \;
        f_n \in \Lip (X) \cap L^2 (\m),
        \mbox{$f_n \to f$ in $L^2 (\m)$}
\right\}.
\]
We say $f \in \mathcal{D} ( \mathsf{Ch} )$
if $f \in L^2 (\mathfrak{m})$ and $\mathsf{Ch} (f) < \infty$.
Note that, for $f \in \mathcal{D} ( \mathsf{Ch} )$,
there exists $| D f | : X \to [0,\infty]$, which is called a minimal
weak upper gradient of $f$. See \cite{AGS2} for a more precise definition
and its equivalence with the minimal relaxed gradient \cite[Theorem~6.2]{AGS2}.
It plays the role of the modulus of gradient of $f$ in the theory of
Sobolev spaces. For instance, it satisfies
\[
\mathsf{Ch} (f) = \frac12 \int_X | D f |^2 \, \d \mathfrak{m}.
\]
We call $(X,d,\m)$ infinitesimally Hilbertian
if $\mathsf{Ch}$ is quadratic form. That is,
$\mathsf{Ch}$ satisfies the parallelogram law (see \cite{AGS3}).
It implies that $f \mapsto | D f |^2$ also becomes an quadratic form.
That is, there exists a bilinear form
$\langle D \cdot, D \cdot \rangle :
\mathcal{D} (\mathsf{Ch}) \times \mathcal{D} (\mathsf{Ch})
\to L^1 (\m)$ such that $\langle D f , D f \rangle = | D f |^2$.
On an infinitesimally Hilbertian $(X,d,\m)$,
we denote the bilinear form corresponding to $2 \mathsf{Ch}$
by $\mathcal{E}$
with $\mathcal{D} (\mathcal{E} ) = \mathcal{D} (\mathsf{Ch})$:
That is,
\[
\mathcal{E} (f,g) = \int_X \langle D f , D g \rangle \, \d \m
\]
and hence $\mathcal{E}(f,f) = 2 \mathsf{Ch} (f)$ and
$\langle D f , D g \rangle$ becomes
the carr\'e du champ associated with $\mathcal{E}$.
We denote the (linear) self-adjoint operator on $L^2 (\m)$
associated with $( \mathcal{E}, \mathcal{D} (\mathcal{E}) )$
by $\Delta$ and the (linear) semigroup of contractions
generated by $\Delta$ by $P_t : L^2 (\m) \to L^2 (\m)$.
Note that we can define $\Delta$ and $P_t$
as a non-linear operator even in absence of
infinitesimal Hilbertianity; see \cite[Section~4]{AGS2}.

We call $(X,d,\m)$ an ${\sf RCD}^*(K,N)$ space for $K \in \R$ and
$N \in (1, \infty]$ if $(X,d,\m)$ is infinitesimally Hilbertian
and $(X,d,\m)$ enjoys the reduced curvature-dimension condition
${\sf CD}^* (K,N)$ introduced in \cite{Bacher:2010bp}.
In \cite{AMS,Erbar:2013wf}, it is shown that, for
infinitesimally Hilbertian $(X,d,\m)$, ${\sf CD}^* (K,N)$ condition
is equivalent to the following three conditions:
\begin{itemize}
\item
There exists $C > 0$ and $x_0 \in X$ such that
\begin{equation} \label{eq:i'ble}
\int_X \e^{ - C d (x_0, x )^2 } \m (\d x) < \infty .
\end{equation}

\item
For $f \in \mathcal{D} (\mathsf{Ch})$
with $| Df | \le 1$ $\m$-a.e.,
$f$ has a 1-Lipschitz representative.

\item
For all $f \in \mathcal{D} (\Delta)$
with $\Delta f \in \mathcal{D} (\mathsf{Ch})$
and $g \in \mathcal{D} ( \Delta ) \cap L^\infty (\m)$
with $g \ge 0$ and $\Delta g \in L^\infty (\m)$,
\begin{multline} \label{eq:Bochner}
\frac12 \int_X | D f |^2 \Delta g \, \d \m
- \int_X \langle D f , D \Delta f \rangle g \, \d \m
\\
\ge
K \int_X | D f |^2 g \, \d \m
+
\frac{1}{N} \int_X ( \Delta f )^2 g \, \d \m.
\end{multline}
\end{itemize}
The last one is nothing but a weak formulation of
Bochner inequality or Bakry-\'Emery's curvature-dimension condition.
See \cite{Ambrosio:2012tp,AGS2,AGS3,AGS_BE-CD} for the case $N = \infty$.
We omit the precise definition of $\mathsf{CD}^* (K,N)$ here
since it is not directly used in this article.
When $K=0$, by definition, $\mathsf{CD}^* (0,N)$
is equivalent to the original curvature-dimension condition
$\mathsf{CD} (0,N)$ in \cite{Sturm_Ric12}. Thus we denote
``$\mathsf{RCD}^*(0,N)$'' by ``$\mathsf{RCD} (0,N)$'' alternatively.
More generally, it is recently proved
in \cite{Cavalletti:2016wn}
that $\mathsf{CD}(K,N)$ is equivalent to $\mathsf{CD}^* (K,N)$
for infinitesimally Hilbertian $(X,d,\m)$
(Indeed, this equivalence is proved under a weaker assumption).
As mentioned in Section~\ref{sec:intro}, a basic
class of $\mathsf{RCD}^* (K,N)$ spaces consists of weighted
Riemannian manifolds satisfying \eqref{eq:wRic}.
It is also known that the $\mathsf{RCD}^* (K,N)$
spaces are stable under the pointed measured Gromov convergence
of metric measure spaces (See \cite{Gigli:2013wi} and references therein;
see Section~\ref{sec:related} also).
In particular,
the limit of a sequence of weighted Riemannian manifolds
satisfying \eqref{eq:wRic} is an $\mathsf{RCD}^* (K,N)$ space.


In the rest of this section, we review several notions and properties
on metric measure spaces, optimal transports and ${\sf RCD}$ spaces
which will be used in the sequel.
We begin with properties on minimal weak upper gradient and the Cheeger energy.
If $f \in \Lip (X) \cap L^2 (\m)$ and $\lip (f) \in L^2 (\mathfrak{m})$,
then $f \in \mathcal{D} ( \mathsf{Ch} )$ and $| D f | \le \lip (f)$.
$f \mapsto | D f |$ is convex in the following sense:
For $f,g \in \mathcal{D} ( \mathsf{Ch} )$ and $\alpha , \beta \in \R$,
\[
| D ( \alpha f + \beta g ) | \le | \alpha | | D f | + | \beta | | D g |.
\]
It indeed implies that $\mathsf{Ch}$ is convex on $L^2 (\mathfrak{m})$.
In addition, $\mathsf{Ch}$ is lower semi-continuous on $L^2 (\mathfrak{m})$
(See \cite[Theorem~4.5]{AGS2}).
If $f$ is constant on a measurable set $A \subset X$, then $| D f | = 0$
$\mathfrak{m}$-a.e.\ on $A$.
Moreover, $| D f | = | D g |$ $\mathfrak{m}$-a.e.\ on $\{ f = g \}$
for $f, g\in \mathcal{D} (\mathsf{Ch})$ (See \cite[Proposition~4.8 (a)(b)]{AGS2}).
We call these properties the locality of the minimal weak upper gradient
in this article.
By using the locality, we can define
$| D f |$ in the extended sense for those measurable $f$ which satisfies
$(-n) \vee ( f \wedge n ) \in \mathcal{D} ( \mathsf{Ch} )$ for each $n > 0$.
See \cite[Section~4]{AGS2}.
Suppose that $(X,d,\mathfrak{m})$ is infinitesimally Hilbertian.
Then we have the Leibniz rule:
For $f, g, h \in \mathcal{D} ( \mathcal{E} ) \cap L^\infty ( \mathfrak{m} )$,
we have $gh \in \mathcal{D} ( \mathcal{E} )$ and
\[
\langle f , g h \rangle
 =
\langle f, g \rangle h
+
\langle f, h \rangle g
\quad \mbox{$\mathfrak{m}$-a.e.}
\]
(See \cite[(4.16)]{Gigli_diff-str} for instance).
Note that
\begin{equation} \label{eq:Dd}
| D d ( x_0 , \cdot ) | = 1
\quad \mbox{$\m$-a.e.}
\end{equation}
holds, where the left hand side is in the extended sense
(see \cite[Proof of Corollary~5.15]{Gigli_diff-str}).
We refer to \cite{AGS2,AGS3,Gigli_diff-str} for other basic properties.

In order to review some properties of the heat flow,
we first recall several notations in optimal transport and metric geometry.
Let $\mathcal{P}_2 (X) \subset \mathcal{P} (X)$ be
the set of probability measure with finite second moment.
That is, $\mu \in \mathcal{P}_2 (X)$ means that
$\| d ( x_0 , \cdot ) \|_{L^2 (\mu)} < \infty$ holds
for some (and hence all) $x_0 \in X$.
For $\mu , \nu \in \mathcal{P} (X)$,
we call $\pi \in \mathcal{P} (X \times X)$ a coupling
of $\mu$ and $\nu$ if the first and second marginal of $\pi$
are $\mu$ and $\nu$ respectively.
That is, for any Borel measurable
$A \subset X$, we have $\pi ( A \times X ) = \mu (A)$
and $\pi ( X \times A ) = \nu (A)$.
We define the $L^2$-Wasserstein distance
$W_2 ( \mu , \nu ) \in [ 0 , \infty ]$ as follows:
\[
W_2 ( \mu , \nu ) : =
\inf
\left\{
\| d \|_{L^2 (\pi)} \mid
\mbox{$\pi$: a coupling of $\mu$ and $\nu$}
\right\}.
\]
Note that $( \mathcal{P}_2 (X) , W_2 )$ becomes
a complete separable geodesic metric space.
Indeed, these properties are inherited from $(X,d,\m)$.
Recall that the convergence in $W_2$ is equivalent
to the weak convergence and the convergence of
the second moment
(See \cite[Theorem~7.12]{book_Vil1} for instance).
We define the relative entropy functional
$\Ent : \mathcal{P}_2 (X) \to ( - \infty , \infty ]$
by
\[
\Ent ( \mu )
 : =
\begin{cases}
\displaystyle \int_X \rho \log \rho \, \d \mathfrak{m}
& \mbox{if $\mu = \rho \mathfrak{m}$},
\\
\infty
& \mbox{if $\mu \not \ll \m$}.
\end{cases}
\]
With the aid of \eqref{eq:i'ble},
$\Ent$ is well-defined as a map as mentioned above.
In addition, $\Ent$ is lower semi-continuous on $( \mathcal{P}_2 (X), W_2 )$
(See \cite[Section~7]{AGS2}).
Let
$
\mathcal{D} ( \Ent )
:=
\{ \mu \in \mathcal{P}_2 (X) \mid \Ent (\mu) < \infty \}
$.
Let $I : \mathcal{P}_2 (X) \to [ 0 , \infty ]$ be the Fisher information
given by
\[
I (\mu) :=
\begin{cases}
\displaystyle 4 \int_X | D \sqrt{\rho} |^2 \, \d \m
&
\mbox{if $\mu = \rho \m , \sqrt{\rho} \in \mathcal{D} ( \mathcal{E} )$},
\\
\infty
&
\mbox{otherwise}.
\end{cases}
\]
Note that we have
\begin{equation} \label{eq:Fisher}
I ( \rho \m ) = \int_X \frac{ |  D \rho |^2 }{ \rho } \, \d \m
\end{equation}
when $I ( \rho \m ) < \infty$ (See \cite[Lemma~4.10]{AGS2}),
where $| D \rho |$ in the right hand side of \eqref{eq:Fisher}
is taken to be an extended sense.
By \cite[Lemma~4.10]{AGS2} again,
$I : \mathcal{P}_2 (X) \to [ 0, \infty ]$ is convex with respect to
convex combinations of elements in $\mathcal{P}_2 (X)$ and
\begin{equation} \label{eq:I-wlsc}
\varliminf_{n \to \infty} I ( \rho_n \mathfrak{m} ) \ge I ( \rho \mathfrak{m} )
\end{equation}
if probability densities $\rho_n$ converges to $\rho$
weakly in $L^1 (\mathfrak{m})$ as $n \to \infty$.
We call a curve $( \gamma_t )_{t \in J}$ indexed by an interval $J \subset \R$
on a metric space $( Y, d_Y )$ absolutely continuous
if there exists $g \in L^1_{\mathrm{loc}} (J)$ such that
\[
d_Y ( \gamma_s , \gamma_t ) \le \int_s^t g (r) \, \d r
\]
for any $s, t \in J$ with $s < t$.
For an absolutely continuous curve $( \gamma_t )_{t \in J}$,
the metric speed $| \dot{\gamma}_t |$
at $t$ is given by
\[
| \dot{\gamma}_t |
: =
\varlimsup_{s \downarrow 0}
\frac{ d_Y ( \gamma_t , \gamma_{t+s} ) }{ s }
\]
Note that we can take $g (s) = | \dot{\gamma}_s |$
in the definition of absolutely continuous curve
if $\gamma$ is absolutely continuous
(See \cite[Theorem~1.1.2]{AGS}).

In the rest of the article,
we always assume $(X,d,\m)$ to be an $\mathsf{RCD} (0,N)$ space
with $N \in [ 1 , \infty )$.
Note that $( \mathcal{E}, \mathcal{D} ( \mathcal{E} ) )$ becomes
a strongly local regular Dirichlet form in the sense of \cite{FOT}
in this framework.
Indeed, it is quasi-regular Dirichlet form by \cite[Theorem~4.1]{Savare:2013wa}
and $X$ is locally compact by the Bishop-Gromov inequality \eqref{eq:BG} below.
The regularity comes from the fact that $\Lip (X) \cap L^2 (\mathfrak{m})$ is
dense in $\mathcal{D} (\mathcal{E})$
(see \cite{AGS_Lip}; see \cite[Proposition~4.10]{AGS3} also).
The chain rule for $\mathcal{E}$ \cite[Theorem~3.2.2]{FOT} says that,
for $f,g \in \mathcal{D} (\mathcal{E}) \cal L^\infty (\mathfrak{m})$
and $\varphi \in C^1 ( \R )$ with $\varphi (0) = 0$, we have
$\varphi (f) \in \mathcal{D} (\mathcal{E}) \cap L^\infty (\mathfrak{m})$ and
\[
\langle \varphi(f) , g \rangle
 =
\varphi' (f) \langle D f, D g \rangle
\quad \mbox{$\mathfrak{m}$-a.e.}
\]

We now turn to review some properties of the heat semigroup $P_t$
which we use in this article.
Since $P_t$ is symmetric and Markovian, there exists an
extension of $P_t$ as an linear contraction from $L^p (\m)$
to itself for $1 \le p \le \infty$ (cf.~\cite[Theorem~4.16]{AGS2}).
In addition, $P_t$ preserves the total mass
(or $( \mathcal{E} , \mathcal{D} (\mathcal{E}) )$ is conservative).
That is, for $f \in L^1 (\m)$ with $f \ge 0$,
$\| P_t f \|_{L^1 (\m)} = \| f \|_{L^1 (\m)}$ holds for $t > 0$.
This is a consequence of \eqref{eq:i'ble} (See \cite[Theorem~4.20]{AGS2}).
Thus $\mu_t = P_t f \m \in \mathcal{P} (X)$ when $f$ is a probability density
with respect to $\mathfrak{m}$ and it can be regarded
as a curve in $\mathcal{P} (X)$ parametrized by $t$.
As a very crucial property of $P_t$ on $\mathsf{RCD}$ spaces,
the curve $( \mu_t )_{t \ge 0}$ as given by $P_t$ in the last sentence
becomes a gradient flow on $\mathcal{P} (X)$
(See \cite{Ambrosio:2012tp,AGS2,AGS3,Erbar:2013wf}).
For any $\mu = \rho \m \in \mathcal{P}_2 (X)$ with $\rho \in L^2 (\m)$,
$\mu_t := P_t \rho \m$ is a gradient flow of $\Ent$
on $( \mathcal{P}_2 (X), W_2 )$ in the sense that $( \mu_t )_{t \ge 0}$
solves $(0,N)$-evolution variational inequality starting from $\mu$
by \cite[Theorem~3.17]{Erbar:2013wf} and its proof.
Indeed, the existence of the solution is one of
characterizing properties of $\mathsf{RCD} (0,N)$ space.
Note that this sort of result is obtained first when $N = \infty$
(See \cite{Ambrosio:2012tp,AGS3}).
For the definition of $(0,N)$-evolution variational inequality,
see \cite[Definition~3.16]{Erbar:2013wf}. We omit the definition
but exhibit some properties obtained from it instead.
First of all, we can extend $P_t$ to be an operator
from $\mathcal{P}_2 (X)$ to itself in the sense that
$P_t \mu = P_t \rho \m$ holds
if $\mu \in \mathcal{P}_2 (X)$ and $\mu = \rho \m$
(See \cite[Theorems~6.1 and 6.2]{Ambrosio:2012tp} for instance).
As an immediate consequence of the definition of
the evolution variational inequality,
$( P_t \mu )_{t > 0}$ is an absolutely continuous curve in
$( \mathcal{P}_2 (X), W_2 )$,
$W_2 ( P_t \mu , \mu ) \to 0$ ($t \downarrow 0$)
and $\Ent ( P_t \mu ) < \infty$
for $\mu \in \mathcal{P}_2 (X)$ and $t > 0$.
By \cite[Remark~3.19]{Erbar:2013wf} and
\cite[Proposition~2.22 (i)]{AGS3},
for $\mu \in \mathcal{P}_2 (X)$,
$t \mapsto \Ent (P_t \mu)$ is absolutely continuous on $( 0, \infty )$
and
$\mu_t = P_t \mu$ solves the energy dissipation identity, i.e.
$\mu_t \to \mu_0$ as $t \to 0$ and for $0 < s < t$,
\begin{equation} \label{eq:E-dissipation0}
\Ent (\mu_s)
=
\Ent ( \mu_t )
+ \frac12 \int_s^t | \dot{\mu}_r |^2 \d r
+ \frac12 \int_s^t I ( \mu_r ) \d r .
\quad \mbox{a.e.~$t$}
\end{equation}
(For instance, see \cite[Definition~2.14]{AGS2} and comments after it).
Since we are on ${\sf RCD} (0,N)$ spaces,
\eqref{eq:E-dissipation0} is equivalent to the following:
\begin{equation} \label{eq:E-dissipation}
- \frac{\d}{\d t} \Ent (\mu_t)
= | \dot{\mu}_t |^2
= I ( \mu_t ) < \infty
\quad \mbox{a.e.~$t$}.
\end{equation}
Here the finiteness follows from \cite[(2.37)]{AGS3}.
Note that \eqref{eq:E-dissipation0} holds
even when $s=0$ if $\mu_0 \in \mathcal{D} ( \Ent )$.

The key property for the proof of the main theorem of this article
is the following space-time $W_2$-control for heat distributions:
For $\mu , \nu \in \mathcal{P}_2 (X)$ and $t, s > 0$,
\begin{equation} \label{eq:control}
W_2 ( P_s \mu , P_t \nu )^2 \le W_2 ( \mu , \nu )^2
+ 2N ( \sqrt{t} - \sqrt{s} )^2 .
\end{equation}
It follows from either the $(0,N)$-evolution variational inequality
(See \cite[Theorem~2.19]{Erbar:2013wf})
or \eqref{eq:Bochner} (See \cite{K14}).

We also review some analytic properties of $P_t$.
As a regularization property,
$P_t f \in \Lip_b (X)$ holds for $f \in L^\infty (\m)$
\cite[Theorem~7.1]{Ambrosio:2012tp}.
One of important tools is
the following Bakry-\'Emery's $L^1$-gradient estimate:
For $f \in \mathcal{D} (\mathcal{E}) \cap L^\infty (\mathfrak{m})$,
we have
\begin{equation} \label{eq:BE1}
| D P_t f | \le \e^{-Kt} P_t ( | D f | ).
\end{equation}
We obtain this bound from the self-improvement of
\eqref{eq:Bochner} with neglecting the term involving $N$
(See \cite{Savare:2013wa} and references therein).
As a technical tool, we recall the following mollification of $P_t$
(See \cite[Section~2.1]{AGS_BE-CD} for instance).
Let $\kappa \in C^\infty_c (( 0, \infty ))$ with $\kappa \geq 0$
and $\int_0^\infty \kappa (r) \, d r = 1$.
For $\eta > 0$ and $f \in L^p (\mu)$ with $p \in [ 1 , \infty ]$,
we define $\mathfrak{h}_\eta f$ by
\[
\mathfrak{h}_{\eta} f
: =
\frac{1}{\eta}  \int_0^\infty
P_r f \; \kappa \left( \frac{r}{\eta} \right)
\, d r.
\]
It is immediate that
$\mathcal{E}( \mathfrak{h}_\eta f - f, \mathfrak{h}_\eta f - f ) \to 0$
and $\| \mathfrak{h}_\eta f - f \|_{L^2(\m)} \to 0$
as $\eta \to 0$ for $f \in \mathcal{D} ( \mathcal{E} )$.
Moreover, for $f \in L^2 ( \mathfrak{m} ) \cap L^\infty ( \mathfrak{m} )$,
$
\mathfrak{h}_\eta f , \Delta (\mathfrak{h}_\eta f)
\in
\mathcal{D} (\Delta) \cap \mathrm{Lip}_b (X)
$.
Here the latter one comes from the following representation:
\[
\Delta \mathfrak{h}_\eta f
=
- \frac{1}{\eta^2} \int_0^\infty
  P_r f \; \kappa' \left( \frac{r}{\eta} \right)
\, d r .
\]
As an additional result from the evolution variational inequality,
$P_t$ admits a symmetric heat kernel $p_t (x,y)$
\cite[Theorem~7.1]{Ambrosio:2012tp}.
That is,
\[
P_t f (x) = \int_X p_t (x,y) f(y) \m ( \d y )
\]
for $f \in L^p (\m)$, $p \in [ 1, \infty ]$.
Moreover, $p_t (x,y)$ admits a sharp two-sided
Gaussian heat kernel estimate \cite{Jiang:2014wo}:
For any $\delta > 0$, there exists $C (\delta) > 0$
such that
\begin{align} \label{eq:HKe}
\frac{1}{C (\delta) \mathfrak{m} ( B_{\sqrt{t}} (x) )}
\exp \left(
- \frac{ d ( x , y )^2 }{( 4 - \delta )t}
\right)
\le
p_t (x,y)
\le
\frac{C (\delta) }{\mathfrak{m} ( B_{\sqrt{t}} (x) )}
\exp \left(
- \frac{ d ( x , y )^2 }{ ( 4 + \delta ) t}
\right)
\end{align}
for $t > 0$ and $x,y \in X$.
Indeed, this is a consequence of the following Li-Yau inequality
\cite[Theorem~1.1]{Jiang:2014vua}:
For $t > 0$ and
$f \in \bigcup_{p \in [ 1, \infty)} L^p (\m)$
with $f \ge 0$ and $f \neq 0$,
\begin{equation} \label{eq:LY}
\frac{ | D P_t f |^2 }{ ( P_t f )^2 }
-
\frac{ \Delta P_t f }{ P_t f }
\le
\frac{N}{2t}
\quad \mbox{$\m$-a.e.}
\end{equation}
%
Note that we can replace $P_t f$ by $p_t (x, \cdot)$, $x \in X$
\cite[Corollary~1.1]{Jiang:2014vua}.

${\sf RCD}^* (K,N)$ spaces satisfy several geometric properties
corresponding to ones on a Riemannian manifold with a lower Ricci curvature bound.
Among them, we recall the Bishop-Gromov inequality on
${\sf RCD} (0,N)$ spaces. It says that, for $0 < r < R$ and $x \in X$,
\begin{equation} \label{eq:BG}
\frac{\m ( B_R (x) )}{ \m ( B_r (x) ) }
\le \left( \frac{R}{r} \right)^N .
\end{equation}
Moreover, \cite{Gig_DePhi} studies the case
when equality holds in \eqref{eq:BG}.
To state their result, we require the following definition:
\begin{Defn}[\protect{\cite[Definition~5.1]{Ketterer:2013tu}}]
A metric measure space $(X' , d' , \m' )$
is a $(0,N-1)$-cone built over a metric measure space $( Y, d_Y , \m_Y )$
if the following holds:
\begin{enumerate}
\item
$X' = [ 0 , \infty ) \times Y / \{ 0 \} \times Y$,

\item
$d' ( [ r , y ], [ s , z ] ) = r^2 + s^2 - 2 rs \cos ( d_Y (y,z) \wedge \pi )$,

\item
$\m' (\d r \d y ) = r^{N-1} \d r \m_Y ( \d y )$.
\end{enumerate}
\end{Defn}

\begin{Thm}[\protect{Special case of \cite[Theorem~1.1]{Gig_DePhi}}]
\label{thm:G-DP}
Let $x \in X$. Suppose that the equality holds in \eqref{eq:BG}
for any $0 < r < R$.
\begin{enumerate}
\item
If $N \ge 2$, then $(X,d,\mathfrak{m})$ is $(0,N-1)$-cone
over an ${\sf RCD}^* ( N-2 , N-1 )$ space
and $x$ is the vertex of the cone.

\item
If $N < 2$, then $(X, d, \mathfrak{m} )$ is isomorphic
to either $( \R, d_{\mathrm{Eucl}} , |x|^{N-1} \d x )$ or
\\
$( [ 0 , \infty ), d_{\mathrm{Eucl}} , x^{N-1} \d x )$,
where $d_{\mathrm{Eucl}}$ is the canonical Euclidean distance.
In both cases, $x \in X$ corresponds to $0$ by the isomorphism.
\end{enumerate}

\end{Thm}
This result will be used in a crucial way
in the proof of our main theorem.

Before closing this section, we re-define of the $\mathcal{W}$-entropy
in a different manner from \eqref{eq:intro-W}.
Indeed we define it as a functional on $\mathcal{P}_2 (X)$.
\begin{Defn}
Let us define the $\mathcal{W}$-entropy
$\mathcal{W}: \mathcal{D} (\Ent) \times ( 0 , \infty ) \to ( - \infty , \infty ]$
by
\[
\mathcal{W} ( \mu , t )
:=
t I (\mu) - \Ent (\mu) - \frac{N}{2} \log t .
\]
\end{Defn}
If $\mu = u \d \mathfrak{m} \in \mathcal{P}_2 (X)$ and
$f$ satisfies $u = \e^{-f} / ( 4 \pi t )^{N/2}$, $\mathcal{W} ( \mu , t )$
coincides with the right hand side of \eqref{eq:intro-W}
up to additive constant (See \cite{LiXD_W-ent0,LiXD_W-ent1,LiXD_W-ent2} for instance).
The choice of additive constant
in \eqref{eq:intro-W} can be regarded as a normalizing constant
so that $\mathcal{W} ( \mu_t , t ) = 0$ on $\R^N$
when $N \in \N$ and $\mu_t$ is a Gaussian kernel.
Since the constant plays no role
in monotonicity and rigidity, we neglect it for brevity.

\section{Monotonicity of $\mathcal{W}$-entropy}
\label{sec:mono}

The goal of this section is to show the monotonicity of the $\mathcal{W}$-entropy
along heat distributions (Theorem~\ref{thm:mono-W}).
For the proof, we show a monotonicity of a rescaled Fisher information
in Proposition~\ref{prop:mono-I2}. Though it is more general than what
we use in the proof of Theorem~\ref{thm:mono-W},
we require this general form in the next section.
We begin with the following auxiliary lemma.
\begin{Lem} \label{lem:mono-I}
For $\mu \in \mathcal{P}_2 (X)$,
$t \mapsto I ( P_t \mu )$ is right continuous and
non-increasing on $[ 0 , \infty )$.
\end{Lem}

\begin{Proof}
We first claim that $t \mapsto I (P_t \mu)$ is non-increasing on $[0,\infty)$.
Let $0 < s < t$ and
we denote the density of $P_s \mu$ with respect to $\mathfrak{m}$
by $\rho$. By \eqref{eq:HKe}, $\rho$ is bounded.
Then, by Bakry-\'Emery's $L^1$-gradient estimate \eqref{eq:BE1}
and the Schwarz inequality for the heat semigroup,
\begin{multline} \label{eq:mono-I0}
I ( P_t \mu )
 =
\int_X
\frac{ | D P_{t-s} \rho |^2 }{ P_{t-s} \rho}
\, \d \mathfrak{m}
 \le
\int_X
\frac{( P_{t - s}( | D \rho | ) )^2}
{P_{t - s} \rho}
\, \d \mathfrak{m}
\\
 \le
\int
P_{t - s} \left(
    \frac{ | D \rho |^2 }{ \rho}
\right)
\, \d \mathfrak{m}
=
\int_X
  \frac{ | D \rho |^2 }{ \rho }
\, \d \mathfrak{m}
=
I ( P_s \mu ).
\end{multline}
Thus the claim holds on $(0, \infty)$.
Let us consider the case $s = 0$.
We may suppose $I (\mu) < \infty$ without loss of generality.
As above, we denote $\mu = \rho \mathfrak{m}$.
For $n \in \N$, let $\bar{\rho}_n := z_n^{-1} \rho_n$,
where $\rho_n := \rho \wedge n$ and
$z_n := \int_X \rho_n \, \d \mathfrak{m}$.
By the locality of the minimal weak upper gradient,
we have $| D \rho_n | = | D \rho | 1_{ \{ \rho < n \} }$ $\mathfrak{m}$-a.e.
Thus we have
\[
I ( \rho \mathfrak{m} )
\ge
\int_X
  \frac{ | D \rho |^2 }{ \rho } 1_{ \{ \rho < n \} }
\, \d \mathfrak{m}
\ge
\frac{1}{n}
\int_X
  | D \rho |^2 1_{ \{ \rho < n \} }
\, \d \mathfrak{m}
=
\frac{1}{n}
\int_X
  | D \rho_n |^2
\, \d \mathfrak{m}
.
\]
Hence $\bar{\rho}_n \in \mathcal{D} ( \mathcal{E} )$ holds
since $I ( \rho \mathfrak{m} ) < \infty$ and $\rho_n \in L^2 (\mathfrak{m})$.
Thus, by the same argument as in \eqref{eq:mono-I0},
\begin{equation} \label{eq:mono-I1}
I ( ( P_t \bar{\rho}_n ) \mathfrak{m} )
\le
\int_X \frac{ | D \bar{\rho}_n |^2 }{ \bar{\rho}_n } \, \d \mathfrak{m}
=
\frac{1}{z_n} \int_X
  \frac{ | D \rho |^2 }{ \rho }
  1_{ \{ \rho < n \} }
\, \d \mathfrak{m}
.
\end{equation}
Since we have \eqref{eq:I-wlsc}, $I ( P_t \mu ) \le I ( \mu )$ holds
by letting $n \to \infty$ in \eqref{eq:mono-I1}.

For a probability density $\rho$ on $X$,
$t \mapsto P_t \rho$ is continuous in $L^1 (\mathfrak{m})$.
Thus \eqref{eq:I-wlsc} yields that
$t \mapsto I (P_t \mu)$ is lower semi-continuous.
It implies the desired right continuity
since $t \mapsto I ( P_t \mu )$ is non-increasing.
\end{Proof}

We next show the following monotonicity formula, which is
closely related with the monotonicity of $\mathcal{W}$-entropy.

\begin{Prop} \label{prop:mono-I2}
\begin{enumerate}
\item
For $t > s \ge 0$, $t' \ge 0$ and $\alpha \in \R$,
\begin{align} \label{eq:pre-mono}
( t + t' )^{2\alpha} I ( P_t \mu )
\le
( s + t' )^{2\alpha} I ( P_s \mu )
+
\frac{N}{2} \cdot
\frac{ ( ( t + t' )^\alpha - ( s + t' )^\alpha )^2 }{ t - s },
\end{align}
where we regard $( s + t' )^{2 \alpha} I ( P_s \mu )$ as 0 if $s = t' = 0$.
\item
$( t + t' )^2 I ( P_t \mu ) - N t /2$
is non-increasing in $t \in ( 0 , \infty )$.
\end{enumerate}
\end{Prop}
\begin{Proof}
(i)
Let $J \subset ( 0 , \infty )$ be the set of $t \in ( 0, \infty )$
satisfying \eqref{eq:E-dissipation}.
For $t,s \in J$ with $t > s$ and $\delta > 0$,
\eqref{eq:control} yields
\begin{multline*}
\frac{ W_2 ( P_t \mu , P_{t + ( t + t' )^\alpha \delta} \mu )^2 }{ \delta^2 }
\le
\frac{ W_2 ( P_s \mu , P_{s + ( s + t' )^\alpha \delta} \mu )^2 }{ \delta^2 }
\\
+ 2N
\left(
    \frac{
      \sqrt{ ( t - s ) + ( ( t + t' )^\alpha - ( s + t' )^\alpha ) \delta }
      -
      \sqrt{ t - s }
    }{ \delta }
\right)^2 .
\end{multline*}
By letting $\delta \to 0$, we obtain \eqref{eq:pre-mono}.
%
Then \eqref{eq:pre-mono} holds for any $0 \le s < t$
since $J$ is dense in $( 0, \infty )$ and
$t \mapsto I (P_t \mu)$ is right continuous by Lemma~\ref{lem:mono-I}.

(ii)
The assertion immediately follows by applying (i) with $\alpha = 1$.
\end{Proof}

Let us turn to consider the monotonicity of the $\mathcal{W}$-entropy.
For later use in Section~\ref{sec:related} (Theorem~\ref{thm:mono-LSI}),
we make our assertion to be slightly stronger than
the usual form, by inserting an additional parameter $t'$.

\begin{Thm}\label{thm:mono-W}
For any $\mu \in \mathcal{P}_2 (X)$ and $t' \ge 0$,
$\mathcal{W} ( P_t \mu , t + t' )$ is non-increasing in $t \in ( 0 , \infty )$.
In addition, the same monotonicity holds for $t \in [ 0 , \infty )$
if $\mu \in \mathcal{D} (\Ent)$ and $t' > 0$.
\end{Thm}

\begin{Proof}
By Proposition~\ref{prop:mono-I2}, for $0 \le s < t$,
\begin{align*}
( t + t' ) I ( P_t \mu ) - ( s + t' ) I ( P_s \mu )
& =
\left( \frac{1}{t + t'} - \frac{1}{s + t'} \right) (s + t' )^2 I ( P_s \mu )
\\
& \qquad +
\frac{1}{t + t'} ( ( t + t' )^2 I ( P_t \mu ) - ( s + t' )^2 I ( P_s \mu )
\\
& \le
\frac{ (s-t) ( s + t' )}{t + t'} I ( P_s \mu )
+ \frac{N}{2 ( t + t' )} (t-s).
\end{align*}
Suppose $s > 0$ if $\mu \notin \mathcal{D} ( \Ent )$ or $t' = 0$.
Suppose $I (\mu) < \infty$ if $s=0$.
Indeed, when $s=0$ and $I ( \mu ) = \infty$, the conclusion obviously holds.
By \eqref{eq:E-dissipation0} and \eqref{eq:E-dissipation}
(see the comment after \eqref{eq:E-dissipation} if $s=0$),
it yields
\begin{align*}
\mathcal{W} & ( P_t \mu , t + t' )
 -
\mathcal{W} ( P_s \mu , s + t' )
\\
& \le
- \frac{ s + t' }{ t + t' } I ( P_s \mu ) (t-s)
+ \frac{N}{2 ( t + t' )}{(t-s)}
- \Ent ( P_t \mu ) + \Ent ( P_s \mu )
- \frac{N}{2} \int_{s + t'}^{t + t'} \frac{\d u}{u}
\\
& =
\int_s^t
\left(
  I ( P_u \mu ) - \frac{s + t'}{t + t'} I ( P_s \mu )
\right)
\, \d u
+ \frac{N}{2}
\int_s^t
\left(
  \frac{1}{t + t'} - \frac{1}{u + t'}
\right)
\, \d u
\\
& \le
\int_s^t \left( I ( P_u \mu ) - I ( P_s \mu ) \right) \d u
+ \frac{ ( t - s )^2 }{t + t'} I ( P_s \mu ).
\end{align*}
Fix $s' \in [ 0 , s ]$ (let $s' > 0$ if $I ( \mu ) = \infty$).
By Lemma~\ref{lem:mono-I},
\begin{align} \label{eq:mono3}
\mathcal{W} ( P_t \mu , t + t' ) - \mathcal{W} ( P_s \mu , s + t' )
\le
\frac{ (t-s)^2 }{ t + t' } I ( P_{s'} \mu ).
\end{align}
Let $n \in \N$ and $t_k = s + k (t-s)/n$ ($k = 0 , 1, \ldots , n$).
By applying \eqref{eq:mono3}
for $( t_k , t_{k-1} )$ instead of $( t, s )$,
we obtain
\begin{align*}
\mathcal{W} ( P_t \mu , t + t' )
 -
\mathcal{W} ( P_s \mu , s + t' )
& =
\sum_{k=1}^n
\mathcal{W} ( P_{t_k} \mu , t_k + t' )
-
\mathcal{W} ( P_{t_{k-1}} \mu , t_{k-1} + t' )
\\
& \le
\frac{( t - s )^2 I ( P_{s'} \mu )}{n ( s + t' )}.
\end{align*}
Hence the conclusion follows by letting $n \to \infty$.
\end{Proof}

\section{Rigidity of $\mathcal{W}$-entropy}
\label{sec:rigid}

Our goal in this section is to show the following theorem.

\begin{Thm} \label{thm:rigid-W}
Suppose that the right upper derivative of $\mathcal{W} ( P_t \mu , t )$ is 0
at $t = t_* \in ( 0 , \infty )$,
that is,
\begin{equation*}
\varlimsup_{t \downarrow t_*}
\frac{ \mathcal{W} ( P_t \mu , t ) - \mathcal{W} ( P_{t_*} \mu , t_* ) }{ t - t_* }
= 0
\end{equation*}
for some $\mu \in \mathcal{P}_2 (X)$ and $t_* \in ( 0 , \infty )$.
Then $\mu = \delta_{x_0}$ for some $x_0 \in X$.
Moreover, we have the following:

\begin{enumerate}
\item
If $N \ge 2$, then $(X,d,\mathfrak{m})$ is $(0,N-1)$-cone
built over an ${\sf RCD}^* ( N-2 , N-1 )$ space
and $x_0$ is the vertex of the cone.

\item
If $N < 2$, then $(X, d, \mathfrak{m} )$ is isomorphic
to either $( \R, d_{\mathrm{Eucl}} , |x|^{N-1} \d x )$ or
\\
$( [ 0 ,n \infty ), d_{\mathrm{Eucl}} , x^{N-1} \d x )$,
where $d_{\mathrm{Eucl}}$ is the canonical Euclidean distance.
In both cases, $x_0 \in X$ corresponds to $0$ by the isomorphism.
\end{enumerate}

In each of these cases,
$\mathcal{W} ( P_t \mu , t )$ is a constant function of $t \in ( 0 , \infty )$.
\end{Thm}
Note that the conclusion in Theorem~\ref{thm:rigid-W} completely corresponds
to the one in Theorem~\ref{thm:G-DP}.
Indeed, we will reduce the proof of Theorem~\ref{thm:rigid-W}
to verification of the assumption of Theorem~\ref{thm:G-DP}.


We begin with the following two lemmas concerning with the Fisher information.

\begin{Lem} \label{lem:less-I}
For any $t > 0$ and $\mu \in \mathcal{P}_2 (X)$,
\begin{equation} \label{eq:less-I}
I ( P_t \mu ) \le \frac{N}{2t}.
\end{equation}
\end{Lem}

This lemma immediately follows
by applying Proposition \ref{prop:mono-I2} (i) with $s=0$ and $\alpha = 1$.
Alternatively, we can show Lemma~\ref{lem:less-I}
by using the Li-Yau inequality \eqref{eq:LY}
(See the proof of Proposition~\ref{prop:id-d}).

\begin{Lem} \label{lem:I-convex}
\begin{enumerate}
\item
$\mu \mapsto I ( P_s \mu )$ is lower semi-continuous
on $\mathcal{P}_2 (X)$ for $s > 0$.

\item
Let $(J , \mathcal{J}, \nu )$ be a probability space.
Let $( \mu_j )_{j \in J} \subset \mathcal{P}_2 (X)$ be a family of probability measures
such that $j \mapsto \mu_j (A)$ is measurable for each measurable $A \subset X$.
Let $\mu_* = \int_J \mu_j \nu ( \d j )$. Suppose $\mu_* \in \mathcal{P}_2 (X)$.
Then, for $s > 0$, we have
\[
I ( P_s \mu_* ) \le \int_J I ( P_s \mu_j ) \nu ( \d j ).
\]
\end{enumerate}
\end{Lem}

\begin{Proof}
(i) By the Lipschitz regularization property of $P_{s/2}$,
$P_{s/2} f \in \Lip_b (X)$ holds for $f \in L^\infty (\mathfrak{m})$.
Take $\mu , \mu^{(n)} \in \mathcal{P}_2 (X)$ ($n \in \N$) and
suppose $W_2 ( \mu^{(n)} , \mu )$ tends to 0 as $n \to \infty$.
By \eqref{eq:control}, we have $W_2 ( P_{s/2} \mu^{(n)} , P_{s/2} \mu ) \to 0$
and thus $P_{s/2} \mu^{(n)}$ converges to $P_{s/2} \mu$ weakly.
Hence, for each $f \in L^\infty (\mathfrak{m})$,
\[
\int_X f \, \d P_s \mu^{(n)}
=
\int_X P_{s/2} f \, \d P_{s/2} \mu^{(n)}
\to
\int_X P_{s/2} f \, \d P_{s/2} \mu
=
\int_X f \, \d P_s \mu .
\]
It implies that the density of $P_s \mu^{(n)}$ converges to that of $P_s \mu$
weakly in $L^1 (\mathfrak{m})$. Then \eqref{eq:I-wlsc} yields the assertion.

(ii) Let $( Z_k )_{k \in \N}$ be $J$-valued, independent and identically distributed
random variables with the law $\nu$.
By the law of large numbers, we have
\[
\lim_{n \to \infty}
W_2 \left(
    \frac{1}{n} \sum_{k=1}^n \mu_{Z_k} , \mu_*
\right)
= 0
\]
almost surely.
Since $I$ is convex in the sense as mentioned in Section~\ref{sec:frame},
the assertion (i) yields
\[
I ( P_s \mu_* )
\le
\varliminf_{n \to \infty}
 I \left( \frac{1}{n} \sum_{k=1}^n P_s \mu_{Z_k} \right)
\le
\varliminf_{n \to \infty}
\frac{1}{n} \sum_{k=1}^n
I \left( P_s \mu_{Z_k} \right)
\]
almost surely.
By taking the expectation and applying the Fatou lemma
in the last inequality,
we obtain the desired inequality.
\end{Proof}
The next lemma shows that the assumption of Theorem~\ref{thm:rigid-W}
implies an identity for the Fisher information.
\begin{Lem} \label{lem:id-I}
Suppose that the assumption of Theorem~\ref{thm:rigid-W} is satisfied.
Then the following holds:
\begin{enumerate}
\item
The equality holds in \eqref{eq:less-I}
for any $t \in ( 0 , t_* ]$.

\item
For $\mu$-a.e.~$x_0 \in X$,
\begin{equation} \label{eq:id-I}
I ( P_t \delta_{x_0} ) = \frac{N}{2t}
\end{equation}
holds for $t \in ( 0, t_* ]$.
In particular, there exists $x_0 \in X$ satisfying \eqref{eq:id-I}.
\end{enumerate}
\end{Lem}

\begin{Proof}
(i) Let $h (r) : = N/2 - r I ( P_r \mu )$.
Note that $h(r) \ge 0$ holds by Lemma~\ref{lem:less-I}.
By the definition of $\mathcal{W}$ and \eqref{eq:E-dissipation},
for $t > t_*$,
\begin{align} \nonumber
\mathcal{W} ( P_{t} \mu , t ) - \mathcal{W} ( P_{t_*} \mu , t_* )
& =
t I ( P_{t} \mu ) - t_* I ( P_{t_*} \mu )
+ \int_{t_*}^{t} \left( I ( P_r \mu ) - \frac{N}{2r} \right) \, \d r
\\ \nonumber
& =
h (t_*) - h (t) - \int_{t_*}^{t} \frac{h(r)}{r} \, \d r
.
\end{align}
Thus the assumption yields
\begin{align}
\label{eq:rigid1}
\varlimsup_{t \downarrow t_*}
\frac{1}{t - t_*} \left(
    h (t_*) - h (t)
    - \int_{t_*}^{t} \frac{h(r)}{r} \, \d r
\right)
= 0 .
\end{align}
By Proposition~\ref{prop:mono-I2} (i), for $\alpha \in ( 0 , 1 )$,
\begin{multline*}
t_*^{ 2\alpha -1 } ( h(t_*) - h(t) )
+ ( t^{2 \alpha -1 } - t_*^{2 \alpha - 1} ) t I (P_t \mu)
\\
=
t^{2 \alpha } I ( P_t \mu ) - t_*^{2\alpha} I ( P_{t_*} \mu )
\le
\frac{N}{2} \cdot \frac{ ( t^\alpha - t_*^\alpha )^2 }{ t - t_* }.
\end{multline*}
Then, combining this with \eqref{eq:rigid1} with Lemma~\ref{lem:mono-I} in mind,
we obtain
\[
t_*^{2\alpha - 2} h ( t_* )
+ ( 2 \alpha -1 ) t_*^{2 \alpha - 1} I ( P_{t_*} \mu )
\le
\frac{N}{2} \cdot \alpha^2 t_*^{2 \alpha - 2} .
\]
By using $I ( P_{t_*} \mu ) = N/2 - h (t_*)$,
this inequality can be simplified to $h (t_*) \le N ( 1 - \alpha )/4$.
Since $\alpha \in ( 0 , 1 )$ is arbitrary,
we obtain $h (t_*) = 0$.
%
%
%
%
By Proposition~\ref{prop:mono-I2} and Lemma~\ref{lem:less-I},
for $0 \le t < t_*$, we have $0 \le t h (t) \le t_* h (t_*)$.
Therefore, $h \equiv 0$ on $( 0, t_* ]$
and this is nothing but the assertion.

(ii)
It suffices to show that $\mu$-a.e. $x_0 \in X$
satisfies $I ( P_{t_*} \delta_{x_0} ) = N / (2t_*)$.
Indeed, the conclusion follows from this by the same argument
as in the proof of the assertion (i).
By Lemmas~\ref{lem:I-convex}~(ii) and \ref{lem:less-I}, we obtain
\[
\frac{N}{2 t_*} = I ( P_{t_*} \mu ) \le \int_X I ( P_{t_*} \delta_x ) \mu ( \d x )
\le
\frac{N}{2t_*}.
\]
Thus all the inequality in the last line must be equality.
It immediately implies the conclusion.
\end{Proof}

For $y \in X$, we define $d_y : X \to [ 0, \infty )$ by $d_y (x) := d (y,x)$.
By employing analytic tools for the heat flow,
the conclusion of the last lemma is transformed into
the crucial identity for the distance function.

\begin{Prop} \label{prop:id-d}
Suppose that there exist $t_* > 0$ and $x_0 \in X$
such that \eqref{eq:id-I} holds for $t \in (0 , t_*]$.
Then,
for $f \in \mathcal{D} (\mathcal{E}) \cap L^1 (\mathfrak{m})$
with $d_{x_0} f , d_{x_0} | D f | \in L^1 ( \mathfrak{m} )$,
we have
\begin{equation} \label{eq:id-d}
- \int_X
\langle
  D d_{x_0}^2 ,
  D f
\rangle
\, \d \mathfrak{m}
=
2N \int_X f \, \d \mathfrak{m}.
\end{equation}
\end{Prop}
\begin{Rem} \label{rem:id-d}
Proposition \ref{prop:id-d} asserts $\Delta d_{x_0}^2 = 2N$
in a distributional sense.
This means that the equality is attained
in the Laplacian comparison theorem on
spaces with $\Ric \ge 0$ and $\dim \le N$.
If we are on a Riemannian manifold,
this identity already implies $X \simeq \R^N$.
\end{Rem}

\begin{Proof}
Take a cut-off function $g_R$ for $R > 1$ by $g_R = 1 \wedge ( R - d_{x_0} )_+ $.
By definition, $0 \le g_R \le 1$, $g_R = 0$ on $B_R (x_0)^c$
and $g_R \to 1$ as $R \to \infty$ pointwisely.
By definition,
$
g_R \in L^1 ( \mathfrak{m} ) \cap L^\infty (\mathfrak{m})
\subset L^2 (\mathfrak{m})
$.
Let $\rho_t$ be the density of $P_t \delta_{x_0}$ with respect to $\mathfrak{m}$
That is, $\rho_t (x) = p_t (x_0 , x )$.
We claim that the equality holds
in the Li-Yau inequality \eqref{eq:LY} for $\rho_t$.
Let $[t , t' ] \subset ( 0 , t_*]$ and set
\[
A_\delta
: =
\left\{
    ( x, s ) \in [ t , t' ] \times X
    \, \left| \,
    \frac{| D \rho_s |^2}{ \rho_t^2 } - \frac{ \Delta \rho_s }{ \rho_s }
    \le \frac{N}{2s} - \delta
    \right.
\right\}
\]
for $\delta > 0$.
By integrating \eqref{eq:LY} for $\rho_s$
by $g_R \d P_s \delta_{x_0} \otimes \d s$ on $X \times [ t, t' ]$,
we obtain
\begin{multline*}
\int_t^{t'} \int_X \frac{ | D \rho_s |^2 }{ \rho_s } g_R \, \d \mathfrak{m} \d s
-
\int_t^{t'} \int_X \Delta \rho_s g_R \, \d \mathfrak{m} \d s
\\
\le
\frac{N}{2} \int_{X \times [t,t'] \setminus A_\delta} \frac{g_R}{s} \, \d P_s \delta_{x_0} \d s
+
\int_{A_\delta} \left( \frac{N}{2s} - \delta \right) g_R \, \d P_s \delta_{x_0} \d s .
\end{multline*}
Since we have
\[
\lim_{R \to \infty} \int_t^{t'} \int_X \Delta \rho_s g_R \, \d \mathfrak{m} \d s
=
\lim_{R \to \infty}
\left(
    \int_X \rho_{t'} g_R \, \d \mathfrak{m}
    -
    \int_X \rho_t g_R \, \d \mathfrak{m}
\right)
= 0 ,
\]
the last inequality implies
\[
\int_t^{t'} I ( \rho_s ) \, \d s
=
\int_t^{t'} \int_X \frac{ | D \rho_s |^2 }{ \rho_s } \, \d \mathfrak{m} \d s
\le
\frac{N}{2} \int_t^{t'} \frac{ \d s }{s}
- \delta
\int_{A_\delta} \, \d P_s \delta_{x_0} \d s
\]
by letting $R \to \infty$.
Since we have \eqref{eq:id-I}, \eqref{eq:HKe} implies that
$A_\delta$ is of null measure with respect to $\d \mathfrak{m} \otimes \d s$.
Hence the equality holds in \eqref{eq:LY} for $\rho_t$ for a.e.~$(x,t)$.
Let $J \subset ( 0 , t_* ]$ be the set of $t$ such that the equality holds
in \eqref{eq:LY} for $\rho_t$ $\mathfrak{m}$-a.e. By the Fubini theorem,
$( 0 , t_* ] \setminus J$ is of null Lebesgue measure.
We show the assertion by combining this identity
with a short time asymptotic of $\rho_t$.
Let $f_0 \in \mathcal{D} ( \Delta ) \cap \mathrm{Lip}_b (X)$
with a bounded support satisfying
$\Delta f_0 \in L^1 (\mathfrak{m}) \cap L^\infty (\mathfrak{m})$.
By \eqref{eq:HKe},
there exists $c > 0$ such that $\rho_t \ge c$ holds on $\supp f_0$.
On the basis of this fact,
the derivation property and the integration by parts formula yields
\begin{align} \nonumber
\int_X \log ( \rho_t ) \Delta f_0 \, \d \mathfrak{m}
& =
- \int_X
\left\langle
  \frac{ D \rho_t}{\rho_t}  ,
  D f_0
\right\rangle
\, \d \mathfrak{m}
\\ \nonumber
& =
-
\int_X \left\langle
  D \rho_t ,
  D \left( \frac{f_0}{ \rho_t } \right)
\right\rangle
\, \d \mathfrak{m}
-
\int_X
\frac{| D \rho_t |^2 }{\rho_t^2} f_0
\, \d \mathfrak{m}
\\ \nonumber
& =
\int_X
\left(
    \frac{ \Delta \rho_t }{ \rho_t }
    -
    \frac{ | D \rho_t |^2 }{ \rho_t^2 }
\right)
 f_0
\, \d \mathfrak{m} .
\end{align}
By plugging the $\mathfrak{m}$-a.e.\ equality in \eqref{eq:LY}
for $t \in J$ into the last identity, we obtain
\begin{align*}
\int_X \log ( \rho_t ) \Delta f_0 \, d \mathfrak{m}
& =
- \frac{N}{2t}
\int_X
 f_0
\, \d \mathfrak{m} .
\end{align*}
By \eqref{eq:HKe}, $\log \rho_t$ is bounded on $\supp \Delta f_0$.
Moreover \eqref{eq:HKe} and \eqref{eq:BG} yield
the Varadhan type short time asymptotic for the heat kernel:
\[
4t \log \rho_t (x) \to - d ( x_0 , x )^2
\quad \mbox{as $t \downarrow 0$ uniformly on each bounded set}.
\]
Thus we obtain
\[
- \frac{N}{2} \int_X f_0 \, \d \mathfrak{m}
=
\lim_{ t \downarrow 0 } \int_X
 ( t \log \rho_t ) \Delta f_0
\, \d \mathfrak{m}
=
- \frac14 \int_X d_{x_0}^2 \Delta f_0 \, \d \mathfrak{m}
=
\frac14 \int_X \langle D d_{x_0}^2 , D f_0 \rangle \, \d \mathfrak{m}.
\]
By combining them, we obtain \eqref{eq:id-d} for $f_0$.

To prove the assertion, we employ an approximation argument.
Let $f \in \mathcal{D} (\mathcal{E}) \cap L^1 (\mathfrak{m})$
and $f_\eta := \mathfrak{h}_\eta ( f \land \eta^{-1} )$
for $\eta > 0$.
By the Lipschitz regularization property of $P_t$,
we have $f_\eta \in \Lip_b (X)$.
By \cite[Lemma~6.7]{Ambrosio:2016ce},
for $R > 0$, there exists $g_R \in \mathrm{Lip}_b (X)$
with $0 \le g_R \le 1$, $g_R|_{B(x_0 , R)} \equiv 1$, $g_R|_{B (x_0 , R+1)} \equiv 0$,
$g_R \in \mathcal{D} ( \Delta )$ and $\Delta g_R \in L^\infty (\mathfrak{m})$.
Note that $g_R f_\eta \in \mathcal{D} (\Delta)$.
Indeed, by applying the Leibniz rule and the integration by parts formula,
for $h \in \mathcal{D} ( \mathcal{E} ) \cap L^\infty (\mathfrak{m})$,
we have
\[
\mathcal{E} ( h , g_R f_\eta )
=
- \int_X
h
\left(
    f_\eta \Delta g_R
    +
    2 \langle D f_\eta , D g_R \rangle
    +
    g_R \Delta f_\eta
\right)
\, \d \mathfrak{m}
\]
(cf.\ \cite[Theorem~4.29]{Gigli_diff-str}).
Thus the same holds for $h \in \mathcal{D} ( \mathcal{E} )$ by a truncation argument.
It implies $g_R f_\eta \in D ( \Delta )$ and
$\Delta ( g_R f_\eta ) = f_\eta \Delta g_R + 2 \langle D f_\eta , D g_R \rangle + g_R \Delta f_\eta$.
From this expression, we can easily verify
$\Delta ( g_R f_\eta ) \in L^1 ( \mathfrak{m}) \cap L^\infty (\mathfrak{m})$.
Thus $g_R f_\eta$ satisfies all assumptions for $f_0$ above
and hence again by the Leibniz rule,
\[
\int_X
\langle
  D d_{x_0}^2 ,
  D g_R
\rangle
f_\eta
\, \d \mathfrak{m}
+
\int_X
\langle
  D d_{x_0}^2 ,
  D f_\eta
\rangle
g_R
\, \d \mathfrak{m}
=
2N \int_X g_R f_\eta \, \d \mathfrak{m}.
\]
By virtue of the fact that $\supp g \subset B_{R+1} (x_0)$,
the last identity implies
\[
\int_X
\langle
  D d_{x_0}^2 ,
  D g_R
\rangle
f
\, \d \mathfrak{m}
+
\int_X
\langle
  D d_{x_0}^2 ,
  D f
\rangle
g_R
\, \d \mathfrak{m}
=
2N \int_X g_R f \, \d \mathfrak{m}
\]
by letting $\eta \to 0$.
Thus, for those $f$ in the assertion,
the conclusion follows by letting $R \to \infty$
in the last identity with the aid of
the locality of minimal weak upper gradient
and \eqref{eq:Dd}.
\end{Proof}

From now on, we study the consequence of the equality \eqref{eq:id-d}.
We set $V(r) := \mathfrak{m} ( B_r (x_0) )$.
We begin with the following two auxiliary lemmas
(Lemmas~\ref{lem:hat-P2} and \ref{lem:vol}).
\begin{Lem} \label{lem:hat-P2}
$\ds
\int_X
 d_{x_0}^p
 \exp \left(
     - c d_{x_0}^2
 \right)
\, \d \mathfrak{m}
< \infty
$ for any $p \ge 0$ and $c > 0$.
\end{Lem}

\begin{Proof}
It suffices to show
$\ds
\exp \left( -  c d_{x_0}^2 \right)
\in L^1 (\mathfrak{m})
$.
By the Fubini theorem,
\begin{equation*}
\int_X
\exp \left(
    - c d_{x_0}^2
\right)
\, \d \mathfrak{m}
 =
\int_X
\left(
    \int_{d_{x_0}}^\infty
        2 c r
    \e^{ - c r^2 }
    \, \d r
\right)
\, \d \mathfrak{m}
 =
2c \int_0^\infty
r V(r)
\e^{ - c r^2 }
\, \d r .
\end{equation*}
By the Bishop-Gromov inequality \eqref{eq:BG},
there exists $C > 0$ such that $V (r) \le V(1) \vee C r^N$.
Thus the desired result holds.
\end{Proof}

\begin{Lem} \label{lem:vol}
Let $f : [ 0 , \infty ) \to \R$ be a measurable function such that
\[
\int_0^\infty r^n |f ( r )| \e^{-r/2} \, \d r < \infty
\]
for all $n \in \N$.
Suppose
\begin{equation} \label{eq:vol-orth}
\frac{N + 2}{2} \int_0^\infty f(r) \e^{- \xi r } \, \d r
= \xi \int_0^\infty r f (r) \e^{- \xi r } \, \d r
\end{equation}
holds for $\xi \in (1/2, 2)$.
Then there exists $c_1 \in \R$ such that
$f(r) = c_1 r^{N/2}$ for a.e.~$r \in [ 0 , \infty )$.
\end{Lem}

\begin{Proof}
Let us choose $c_1 \in \R$ so that
\[
\int_0^\infty f (r) \e^{-r} \, \d r
= c_1 \int_0^\infty r^{N/2} \e^{-r} \, \d r
\]
and set $g(r): = f(r) - c_1 r^N$.
We can easily verify that \eqref{eq:vol-orth} holds
for $g$ instead of $f$.
That is,
\begin{equation} \label{eq:orth}
\frac{N + 2}{2} \int_0^\infty g(r) \e^{- \xi r } \, \d r
= \xi \int_0^\infty r g (r) \e^{- \xi r } \, \d r .
\end{equation}
We will show
\begin{equation} \label{eq:orth1}
\int_0^\infty r^n g(r) \e^{-r} \, \d r = 0
\end{equation}
by induction for $n \in \N \cup \{ 0 \}$.
The assertion for $n = 0$ holds by the definition of $g$.
Suppose the claim to be true for $n \in \N$.
By differentiating \eqref{eq:orth} $n$ times
with respect to $\xi$ at $\xi = 1$
with taking the assumption of $f$ into account.
Then, by the assumption of the induction,
we immediately obtain \eqref{eq:orth1} for $n+1$.
Thus \eqref{eq:orth1} holds for any $n \in \N \cup \{ 0 \}$.
It means that $g(r)$ is orthogonal to any polynomial
in $L^2 ( \e^{-r} \d r )$.
Hence $g = 0$ a.e.\ and the conclusion follows.
\end{Proof}

With keeping Lemma~\ref{lem:hat-P2} in mind,
let us define $\hat{\rho}^{x_0} : ( 0 , \infty ) \times X \to \R$ by
\[
\hat{\rho}^{x_0}_t (x)
: =
\frac{1}{Z(t)}
\exp \left(
    - \frac{ d_{x_0} (x)^2 }{4t}
\right),
\]
where $Z(t)$ is a normalizing constant so that
$\| \hat{\rho}^{x_0} \|_{L^1 (\m)} = 1$. For simplicity of notations
we denote $\hat{\rho}^{x_0}$ by $\hat{\rho}$ if there is no
possibility of confusions.
Note that $\hat{\mu}_t : = \hat{\rho}_t \mathfrak{m}
\in \mathcal{P}_2 (X)$ holds by Lemma~\ref{lem:hat-P2}.

\begin{Lem} \label{lem:hat-F}
Suppose \eqref{eq:id-d}.
Then $I ( \hat{\mu}_t ) = N/(2t)$ for $t > 0$ and
there exists $c_* , c_{**} > 0$ such that
$Z(t) = c_* t^{N/2} = c_{**} V (\sqrt{t})$ for $t > 0$.
In particular, $X$ is non-compact.
\end{Lem}

\begin{Proof}
Since $\hat{\rho}_t \in \Lip (X)$, we have
$
\lip( \hat{\rho}_t )
\le
d_{x_0} \hat{\rho}_t /(2t)
$.
Since $d_{x_0} \hat{\rho}_t , d_{x_0}^2 \hat{\rho}_t \in L^2 (\mathfrak{m})$
by Lemma~\ref{lem:hat-P2},
$\hat{\mu}_t \in \mathcal{D} ( I )$ and
\begin{equation} \label{eq:D-hat}
| D \hat{\rho}_t |
\le
\frac{ d_{x_0} \hat{\rho}_t }{2t}
\end{equation}
holds.
For $R > 0$, take $g_R\in \Lip_b (X)$ satisfying
$0 \le g_R \le 1$, $g_R|_{B(x_0 , R)} \equiv 1$, $g_R|_{B (x_0 , R+1)} \equiv 0$.
Then, by the Leibniz rule, the locality of minimal weak upper gradient
and the chain rule, we have
\begin{align} \nonumber
\int_X
\bigg(
  \frac{ | D \hat{\rho}_t |^2  g_R }{ \hat{\rho}_t }
&  +
  \langle D \hat{\rho}_t , D g_R \rangle
\bigg)
\, \d \mathfrak{m}
 =
\int_X
\frac{ \langle D \hat{\rho}_t , D ( \hat{\rho}_t g_R ) \rangle }{ \hat{\rho}_t }
\, \d \mathfrak{m}
\\ \nonumber
& =
-\frac{1}{4t}
\int_X
\langle D d_{x_0}^2 , D ( \hat{\rho}_t g_R ) \rangle
\, \d \mathfrak{m}
\\ \label{eq:pre-id-F-hat}
& =
-\frac{1}{4t}
\int_X
\left(
    \langle D d_{x_0}^2 , D \hat{\rho}_t \rangle g_R
    +
    \langle D d_{x_0}^2 , D g_R \rangle \hat{\rho}_t
\right)
\, \d \mathfrak{m}.
\end{align}
By virtue of \eqref{eq:Dd}, \eqref{eq:D-hat}, Lemma~\ref{lem:hat-P2}
and the locality of minimal weak upper gradient,
letting $R \to \infty$ in \eqref{eq:pre-id-F-hat} implies
\begin{equation} \label{eq:id-F-hat}
I ( \hat{\mu}_t )
=
\int_X
\frac{ | D \hat{\rho}_t |^2 }{\hat{\rho}_t}
\, \d \mathfrak{m}
 =
- \frac{1}{4t}
\int_X
\left\langle
    D d_{x_0}^2
    D \hat{\rho}_t
\right\rangle
\, \d \mathfrak{m}
=
\frac{N}{2t}
\int_X \hat{\rho}_t \, \d \mathfrak{m}
= \frac{N}{2t},
\end{equation}
where the third inequality comes from Proposition~\ref{prop:id-d}.
Thus the first assertion holds.
For the second assertion, we compute
$I ( \hat{\mu}_t )$ in a different manner.
By the Leibniz rule, the chain rule and \eqref{eq:Dd},
\begin{align*}
\int_X
\langle D d_{x_0}^2 , D ( \hat{\rho}_t g_R ) \rangle
\, \d \mathfrak{m}
& =
\int_X
\left(
    - \frac{1}{4t} | D d_{x_0}^2 |^2  \hat{\rho}_t g_R
    +
    \langle D d_{x_0}^2 , D g_R \rangle \hat{\rho}_t
\right)
\, \d \mathfrak{m}
\\
& =
\int_X
\left(
    - \frac{1}{t} d_{x_0}^2 \hat{\rho}_t g_R
    +
    \langle D d_{x_0}^2 , D g_R \rangle \hat{\rho}_t
\right)
\, \d \mathfrak{m}.
\end{align*}
By substituting this identity into \eqref{eq:pre-id-F-hat}
and letting $R \to \infty$, we obtain
\begin{equation} \label{eq:id-F-hat2}
I ( \hat{\mu}_t )
 =
\frac{1}{4t^2}
\int_X
d_{x_0}^2
\hat{\rho}_t
\, \d \mathfrak{m}.
\end{equation}
As in the proof of Lemma~\ref{lem:hat-P2},
we have
\[
\int_X d_{x_0}^{2p}
\exp \left(
    - \frac{ d_{x_0}^2 }{ 4 t }
\right)
\, \d \mathfrak{m}
=
\int_0^\infty
V(r) \left(
    \frac{r^{2p+1}}{2t} - 2p r^{2p-1}
\right)
\exp \left(
    - \frac{ r^2 }{ 4 t }
\right)
\, \d r
\]
for $p \ge 0$. By combining this
with $\hat{\mu}_t (X) = 1$, \eqref{eq:id-F-hat} and \eqref{eq:id-F-hat2},
we obtain
\begin{equation} \label{eq:vol1}
Z(t)
 =
\frac{1}{2t} \int_0^\infty
r V(r) \exp \left( - \frac{r^2}{4t} \right)
\, \d r
 =
\frac{1}{2Nt}
\int_0^\infty
V(r) \left( \frac{r^3}{2t} - 2 r \right)
\exp \left( - \frac{r^2}{4t} \right)
\, \d r  .
\end{equation}
Thus we have
\begin{equation*}
\frac{N+2}{2} \int_0^\infty
r V(r) \exp \left( - \frac{r^2}{4t} \right)
\, \d r
 =
\int_0^\infty
\frac{r^3 V(r) }{4t}
\exp \left( - \frac{r^2}{4t} \right)
\, \d r .
\end{equation*}
After the change of variable $\tilde{r} = r^2$,
we can apply Lemma~\ref{lem:vol} to conclude that
there exists $c_1 \in \R$ such that
$V(r) = c_1 r^N$ for a.e.~$r$.
Note that $c_1 > 0$ holds by the definition of $V$.
Since $V$ is left-continuous, The last identity for $V$
holds for all $r \ge 0$.
Then the assertion for $Z(t)$ follows
from the first identity in \eqref{eq:vol1}.
Finally, $X$ must be non-compact
since $V (r) \to \infty$ as $r \to \infty$.
\end{Proof}

As we see in Theorem~\ref{thm:G-DP}, the conclusion of Lemma~\ref{lem:hat-F}
is already sufficient to specify $(X,d,\m)$ as in Theorem~\ref{thm:rigid-W}
by Lemma~\ref{lem:id-I} and Proposition~\ref{prop:id-d}.
To show $\mu$ to be a Dirac measure, we need the following additional arguments.

\begin{Prop} \label{prop:hat-GF}
Suppose \eqref{eq:id-d}.
Then
$\hat{\mu}_t = P_t \delta_{x_0}$ for $t > 0$.
In particular, $\hat{\mu}_t \to \delta_{x_0}$ as $t \to 0$.
\end{Prop}

\begin{Proof}
Let $\delta > 0$ and $\hat{\mu}^\delta_t := \hat{\mu}_{t + \delta}$.
Note that $\hat{\mu}^\delta_0 = \hat{\mu}_\delta \in \mathcal{D} ( \Ent )$
and $\hat{\rho}_\delta \in L^2 (\mathfrak{m})$.
We first show $\hat{\mu}^\delta_t = P_t \hat{\mu}_\delta$.
By \cite[Theorem~9.3]{AGS2},
$( P_t \hat{\mu}_\delta )_{t \ge 0}$ is a unique gradient flow of $\Ent$
starting from $\hat{\mu}_\delta$
in the sense of the energy dissipation identity \eqref{eq:E-dissipation0}.
Thus it suffices to show that
$( \hat{\mu}^\delta_t )_{ t \ge 0 }$ is also a gradient flow
of $\Ent$ in the same sense.
We show it by following a strategy in \cite{AGS2}.

Since we are on ${\sf RCD} (0,N)$ space,
we already know that the descending slope $| D^- \Ent |$ of $\Ent$
is an upper gradient of $\Ent$ and $| D^- \Ent |^2 = I$
by \cite[Theorem 9.3]{AGS2}.
In this case, the inequality ``$\le$'' in \eqref{eq:E-dissipation0} holds automatically
if $( \mu_r )_{r \ge 0}$ is absolutely continuous curve in $( \mathcal{P}_2 (X), W_2 )$.
Thus the proof is reduced to show the absolute continuity of $( \hat{\mu}^\delta_t )_{t > 0}$
and the following:
\begin{align} \label{eq:Ent-D}
\frac{ \d}{ \d t} \Ent ( \hat{\mu}^\delta_t )
& =
- I ( \hat{\mu}^\delta_t )
\quad \mbox{for a.e.~$t$},
\\
\label{eq:W-speed}
\varlimsup_{s \downarrow 0} \frac{ W_2 ( \hat{\mu}^\delta_t , \hat{\mu}^{\delta}_{t+s} )^2 }{ s^2 }
& \le
I ( \hat{\mu}^\delta_t )
\quad \mbox{for a.e.~$t$}.
\end{align}
By the definition of $\hat{\mu}^\delta_t$,
\eqref{eq:id-F-hat2}, \eqref{eq:id-F-hat} and Lemma~\ref{lem:hat-F},
\begin{align*}
\Ent ( \hat{\mu}^\delta_t )
=
- \frac{1}{4(t+\delta)} \int_X d_{x_0}^2 \hat{\rho}_{t+\delta} \, \d \mathfrak{m}
- \log Z(t+\delta)
=
- \frac{N}{2}
- \log c_* - \frac{N}{2} \log (t+\delta).
\end{align*}
Thus, again by Lemma~\ref{lem:hat-F},
\begin{align*}
\frac{ \d}{ \d t} \Ent ( \hat{\mu}^\delta_t )
& =
- \frac{N}{2 ( t + \delta )}
=
- I ( \hat{\mu}^\delta_t ).
\end{align*}
Hence we obtain \eqref{eq:Ent-D}.
For \eqref{eq:W-speed}, the Kantorovich duality yields
\begin{align} \label{eq:Kant}
\frac{ W_2 ( \hat{\mu}^\delta_t , \hat{\mu}^{\delta}_{t+s} )^2 }{ s^2 }
=
\frac{2}{s}
\sup_{\varphi \in \Lip_b (X)}
\left[
\int_X Q_s \varphi \, \d \hat{\mu}^{\delta}_{t+s}
-
\int_X \varphi \, \d \hat{\mu}^{\delta}_t
\right],
\end{align}
where $Q_s$ is the Hopf-Lax semigroup defined by
\[
Q_s \varphi (x)
: =
\inf_{y \in X}
\left[
\varphi (y)
+ \frac{ d(x,y)^2 }{ 2s }
\right].
\]
It is known that $Q_s \varphi \in \Lip_b (X)$ for $\varphi \in \Lip_b (X)$
(see \cite{BEHM_HJ} for instance).
It is not difficult to verify that we can suppose $\varphi$ to be
of bounded support in the range of the supremum in \eqref{eq:Kant}.
Let $\varphi \in \Lip_b (X)$ with a bounded support.
It is easy to see $\varphi \in \mathcal{D} ( \mathcal{E} )$.
Note that $\hat{\rho}_{t + \delta}$ is differentiable in $t \in ( 0 , \infty )$ and
possesses a sufficiently good integrability by Lemma~\ref{lem:hat-P2}.
Thus, by the dominated convergence theorem,
$\int_X Q_r \varphi \hat{\rho}_{t + \delta + r } \, \d \mathfrak{m}$ is
a.e.\ differentiable in $r$ and
\begin{align} \label{eq:W-speed0}
\int_X Q_s \varphi \, \d \hat{\mu}^\delta_{t+s}
-
\int_X \varphi \, \d \hat{\mu}^\delta_t
=
\int_0^s \left(
    \frac{\d}{\d r} \int_X
    Q_r \varphi \hat{\rho}_{\delta + t + r }
    \, \d \mathfrak{m}
\right)
\, \d r .
\end{align}
Since $Q_r \varphi$ satisfies the following Hamilton-Jacobi equation
\[
\frac{ \partial}{ \partial s } Q_s \varphi + \frac12 \lip( Q_s \varphi )^2 = 0
\]
(see \cite[Theorem~3.6]{AGS2}),
by Lemma~\ref{lem:hat-F} and the fact $|D f| \le \lip (f)$,
we have
\begin{align} \nn
\frac{\d}{\d r} \int_X
Q_r \varphi \hat{\rho}_{\delta + t + r }
\, \d \mathfrak{m}
& \le
- \frac12
\int_X
 | D Q_r \varphi |^2 \hat{\rho}_{t + \delta + r}
\, \d \mathfrak{m}
+
\frac{1}{4 (t + \delta + r )^2}
\int_Z
 d_{x_0}^2 Q_r \varphi \hat{\rho}_{t + \delta + r }
\, \d \mathfrak{m}
\\
\label{eq:W-speed1}
& \hspace{3em}
- \frac{N}{2( t + \delta + r )}
\int_X
 Q_r \varphi \hat{\rho}_{t + \delta + r }
\, \d \mathfrak{m}.
\end{align}
By \eqref{eq:id-d}, the Leibniz rule, the chain rule and \eqref{eq:Dd}
together with a localization argument
as we did in the proof of Lemma~\ref{lem:hat-F},
\begin{align*}
- & \frac{N}{2( t + \delta + r )}
\int_X
 Q_r \varphi \hat{\rho}_{t + \delta + r }
\, \d \mathfrak{m}
 =
\frac{1}{4( t + \delta + r )}
\int_X \left\langle
    D d_{x_0}^2 ,
 D ( Q_r \varphi \hat{\rho}_{t + \delta + r } )
\right\rangle
\, \d \mathfrak{m}
\\
& =
\frac{1}{4 ( t + \delta + r )}
\left(
\int_X \left\langle
    D d_{x_0}^2 ,
 D  Q_r \varphi
\right\rangle
\hat{\rho}_{t + \delta + r }
\, \d \mathfrak{m}
+
\int_X
Q_r \varphi
\left\langle
    D d_{x_0}^2 ,
 D \hat{\rho}_{t + \delta + r }
\right\rangle
\, \d \mathfrak{m}
\right)
\\
& =
-
\int_X \left\langle
    D \hat{\rho}_{t + \delta + r } ,
 D  Q_r \varphi
\right\rangle
\, \d \mathfrak{m}
-
\frac{1}{4 ( t + \delta + r )^2}
\int_X
    d_{x_0}^2 Q_r \varphi
    \hat{\rho}_{t + \delta + r }
\, \d \mathfrak{m} .
\end{align*}
Thus, by substituting this identity
into \eqref{eq:W-speed1}, we obtain
\begin{align*}
\frac{\d}{\d r} \int_X
Q_r \varphi \hat{\rho}_{\delta + t + r }
\, \d \mathfrak{m}
& \le
- \frac12
\int_X
 | D Q_r \varphi |^2
\, \d \hat{\mu}^\delta_{t+r}
-
\int_X
\frac{1}{\hat{\rho}_{t+\delta+r}}
\left\langle
  D \hat{\rho}_{t + \delta + r } ,
  D Q_r \varphi
\right\rangle
\, \d \hat{\mu}^\delta_{t+r}
\\
& \le
\frac12 I ( \hat{\mu}^\delta_{t+r} ) = \frac{N}{4 ( t + \delta + r )}.
\end{align*}
By applying this inequality to \eqref{eq:W-speed0},
\eqref{eq:Kant} yields
\[
\frac{ W_2 ( \hat{\mu}^\delta_t , \hat{\mu}^{\delta}_{t+s} )^2 }{ s^2 }
\le
\frac{1}{s} \int_0^s
\frac{N}{ ( t + \delta + r )}
\, \d r .
\]
This inequality together with Lemma~\ref{lem:hat-F}
yields \eqref{eq:W-speed}.
It is immediate that the last inequality implies
the absolute continuity of $( \hat{\mu}^\delta_t )_{t \ge 0}$.
Hence $\hat{\mu}_{\delta + t} = P_t \hat{\mu}_\delta$ for $t \ge 0$.

Note that we have
\[
W_2 ( \hat{\mu}_\delta , \delta_{x_0} )^2
 =
\int_X d_{x_0}^2 \hat{\rho}_t \, \d \mathfrak{m}.
\]
Here the right hand side goes to 0 as $\delta \to 0$
by virtue of the explicit expression of $Z(t)$ in Lemma~\ref{lem:hat-F}.
This fact implies the last assertion.
In addition, by \eqref{eq:control},
\[
W_2 ( \hat{\mu}_{ t + \delta} , P_t \delta_{x_0} )
=
W_2 ( P_t \hat{\mu}_{\delta} , P_t \delta_{x_0} )
\le
W_2 ( \hat{\mu}_\delta , \delta_{x_0} ).
\]
Then the assertion for $t > 0$ holds by letting $\delta \downarrow 0$
since $( \hat{\mu}_{t} )_{t > 0}$ is a continuous curve
in $( \mathcal{P}_2 (X) , W_2 )$.
\end{Proof}

\begin{Rem}
A converse of Proposition~\ref{prop:hat-GF} holds in the following sense.
Let $\hat{\mu}_t$ be as above
and suppose $Z(t) = c_* t^{N/2}$ for some constant $c_* > 0$.
If $\hat{\mu}_t$ is a solution to the heat equation,
then $\mathcal{W} ( \hat{\mu}_t , t )$ is constant in $t$
and hence $\mathcal{W} ( \hat{\mu}_t , t )$ has
vanishing $t$-derivatives.
\end{Rem}

\begin{Lem} \label{lem:Dirac}
Suppose that the equality holds in \eqref{eq:less-I}.
Then $\mu$ is a Dirac measure.
\end{Lem}

\begin{Proof}
Let $\Phi : ( 0, \infty ) \times [ 0 , \infty )$ be given by
$\Phi ( u, v ) := v^2 / u$. It is verified in a straightforward way
that $\Phi$ is convex and that
\[
\Phi (
  ( 1 - \lambda ) ( u_1 , v_1 )
  +
  \lambda ( u_2 , v_2 )
)
=
( 1 - \lambda ) \Phi ( u_1 , v_1 )
  +
\lambda \Phi ( u_2 , v_2 )
\]
holds for some $\lambda \in (0,1)$ if and only if $v_1 / u_1 = v_2 / u_2$.
Note that we have
\[
P_t \mu
=
\frac12 \left(
    \int_X P_t \delta_x \mu (\d x)
    +
    \int_X P_t \delta_y \mu (\d y)
\right)
=
\int_{X \times X}
  \frac12 ( P_t \delta_x + P_t \delta_y )
\mu^{\otimes 2} (\d x \d y).
\]
By combining Lemma~\ref{lem:I-convex}~(ii) with the convexity of minimal weak upper gradient
and the convexity of $\Phi$, we have
\begin{align*}
\frac{N}{2t} = I ( P_t \mu )
& \le
\int_X
\left(
    \int_{X^2}
    \Phi \left(
      \frac{
          P_t \delta_x
          +
          P_t \delta_y
      }{2},
      \frac{
          | D P_t \delta_x |
          +
          | D P_t \delta_y |
      }{2}
   \right)
   \, \mu^{\otimes 2} (\d x \d y)
\right)
\d \mathfrak{m}
\\
& \le
\int_X
\left(
    \int_{X}
    \Phi \left(
      P_t \delta_x ,
      | D P_t \delta_x |
    \right)
    \, \mu (\d x)
\right)
\d \mathfrak{m}
\\
& =
\int_X
  I ( P_t \delta_x )
\, \mu ( \d x )
\le \frac{N}{2t},
\end{align*}
where the last inequality follows from Lemma~\ref{lem:less-I}.
Hence all the inequalities are indeed equality.
By the property of $\Phi$ mentioned at the beginning of the proof
together with the Fubini theorem,
for $\mu^{\otimes 2}$-a.e.~$(x,y)$,
we have
\begin{align} \label{eq:Dirac1}
\frac{ | D P_t \delta_x | }{ P_t \delta_x }
=
\frac{ | D P_t \delta_y | }{ P_t \delta_y }
\quad
\mbox{$\mathfrak{m}$-a.e.}
\end{align}
By virtue of \eqref{eq:id-I},
Propositions~\ref{prop:id-d} and \ref{prop:hat-GF} ensure
$P_t \delta_x = \hat{\rho}^x_t \mathfrak{m}$ $\mu$-a.e.~$x$.
This representation of $P_t \delta_x$ together with \eqref{eq:Dd}
implies
\begin{align} \label{eq:Dirac2}
\frac{ | D P_t \delta_x | }{ P_t \delta_x }
=
\frac{ | D \hat{\rho}^x_t | }{  \hat{\rho}^x_t }
=
\frac{d_x}{2t}
\quad
\mbox{$\mathfrak{m}$-a.e.}
\end{align}
By combining \eqref{eq:Dirac1} and \eqref{eq:Dirac2},
we conclude that $d_x = d_y$ $\mathfrak{m}$-a.e. for $\mu^{\otimes 2}$-a.e.~$(x,y)$.
Hence $\mu$ must be a Dirac measure.
\end{Proof}

\begin{tProof}{Theorem~\ref{thm:rigid-W}}
By Lemma~\ref{lem:Dirac}, we know $\mu = \delta_{x_0}$
for some $x_0 \in X$. Then, by Lemma~\ref{lem:hat-F},
there exists $c > 0$ such that $\m ( B_r (x_0) ) = c r^N$
holds for all $r > 0$. It yields the equality in \eqref{eq:BG}.
Then the conclusion follows from Theorem~\ref{thm:G-DP}.
\end{tProof}

\section{Related results}
\label{sec:related}

Here we gather some results related with our main theorem.
First we show that the heat flow coincides with the $W_2$-geodesic
in the rigidity case. According to Proposition~\ref{prop:hat-GF},
we define $\hat{\mu}_0 : = \delta_{x_0}$.

\begin{Prop} \label{prop:hat-G}
Suppose \eqref{eq:id-d}.
Then
$( \hat{\mu}_{t^2/(2N)} )_{t \ge 0}$ is a unit-speed minimal geodesic in $W_2$.
\end{Prop}

\begin{Proof}
Let $\bar{\mu}_t := \hat{\mu}_{t^2/(2N)}$.
By \eqref{eq:id-F-hat2} and Lemma~\ref{lem:hat-F},
\begin{align*}
W_2 ( \bar{\mu}_0 , \bar{\mu}_{t} )
=
\left( \int_X d_{x_0}^2 d \bar{\mu}_{t} \right)^{1/2}
=
t.
\end{align*}
On the other hand,
by Proposition~\ref{prop:hat-GF} and Lemma~\ref{lem:hat-F},
\begin{align*}
\varlimsup_{s \downarrow 0}
\frac{W_2 ( \bar{\mu}_t , \bar{\mu}_{t + s} )}{s}
=
\frac{t}{N} \sqrt{ I ( \bar{\mu}_t ) } = 1
\end{align*}
(See \cite[Definition~2.14]{AGS2}).
We denote the left hand side of the last identity
by $| \dot{\bar{\mu}}_t |$.
This is the metric speed for $\bar{\mu}_t$ in $( \mathcal{P}_2 (X) , W_2 )$.
Thus, as remarked in Section~\ref{sec:frame},
we have
\[
W_2 ( \bar{\mu}_s , \bar{\mu}_t ) \le \int_s^t | \dot{\bar{\mu}}_r | \, \d r
\]
for $0 < s < t$ and hence
\begin{align*}
t =
W_2 ( \bar{\mu}_0 , \bar{\mu}_t )
\le
W_2 ( \bar{\mu}_0 , \bar{\mu}_s )
 +
W_2 ( \bar{\mu}_s , \bar{\mu}_t )
\le
s + \int_s^t | \dot{\bar{\mu}}_r | \, \d r
= t .
\end{align*}
Thus all the last inequality must be equality.
In particular, $W_2 ( \bar{\mu}_s , \bar{\mu}_t ) = t - s$.
Therefore the assertion holds.
\end{Proof}

The second result asserts monotonicity in time of the infimum of $\mathcal{W}$-entropy.
It partially extends \cite[Theorem~2.4]{LiXD_W-ent0} to our framework.
It is related with the (logarithmic) Sobolev inequality.
See Remarks~\ref{rem:LSI} and \ref{rem:SI} below.

\begin{Thm} \label{thm:mono-LSI}
Let
$
c(t) :=
\inf \{ \mathcal{W} ( \mu , t ) \mid \mu \in \mathcal{D} ( \Ent ) \}
$
for $t \in ( 0, \infty )$.
Then $c$ is non-increasing.
\end{Thm}

\begin{Proof}
Let $0 < s < t$.
Note that $P_r \mu \in \mathcal{D} ( \Ent )$ holds
for each $\mu \in \mathcal{P}_2 (X)$ and $r > 0$.
Thus, for any $\mu \in \mathcal{D} ( \Ent )$,
we have
\[
\mathcal{W} ( \mu , s )
\ge
\mathcal{W} ( P_{t-s} \mu , t )
\ge c (t)
\]
by Theorem~\ref{thm:mono-W}.
Hence the conclusion holds
by taking infimum over all $\mu \in \mathcal{D} ( \Ent )$
in the left hand side of the last inequality.
\end{Proof}

\begin{Rem} \label{rem:LSI}
By the definition of $c(t)$ in Theorem~\ref{thm:mono-LSI}, we have
\[
\Ent (\mu) \le t I (\mu) - \frac{N}{2} \log t - c(t)
\]
for all $\mu \in \mathcal{D} ( \Ent )$.
That is, if $c (t) > - \infty$ for some $t$,
then we have a defective log-Sobolev inequality.
\end{Rem}

\begin{Rem} \label{rem:SI}
Let $N > 2$.
Then $c: = \inf_{t > 0} c(t) > - \infty$
if and only if there exists a constant
$C > 0$ depending on $c$ and $N$ such that
\[
\exp \left( \frac{2}{N} \Ent ( \mu ) \right )
\le
C I ( \mu )
\]
holds for any $\mu \in \mathcal{D} ( \Ent )$.
By \cite[Theorem 6.2.3]{BGL}, this inequality implies
the Sobolev inequality:
\[
\| f \|_{2N/(N-2)}^2 \le C \int_X | D f |^2 \, \d \mathfrak{m}
\]
Indeed, weaker forms of these inequalities are equivalent.
Note that $c = \lim_{t \to \infty} c (t)$ by Theorem~\ref{thm:mono-LSI}.
Moreover, by \cite[Theorem~6.3.1]{BGL},
A weaker form of the Sobolev inequality is
equivalent to the ultracontractivity:
\[
\| P_t \|_{1 \to \infty} \le \frac{C'}{t^{N/2}}.
\]
Combined with the heat kernel lower bound \eqref{eq:HKe},
$c > - \infty$ implies that there exists $C'' > 0$ such that
$\mathfrak{m} ( B_r (x) ) \ge C'' r^N$ for each $r > 0$ and $x \in X$.
Note that the Bishop-Gromov inequality implies a similar bound but
it is local in the sense it holds for $r < R$ for each fixed $R$,
and the constant corresponding to $C''$ depends on $R$ and $x$.
Note that a similar result is obtained
in \cite[Theorem~6.1]{LiXD_W-ent0}.
It is shown on weighted Riemannian manifolds
but the same argument works even in our framework.
\end{Rem}

The third result asserts a stronger rigidity under a stronger assumption.

\begin{Thm} \label{thm:s-rigid-W}
Suppose that the assumption of Theorem~\ref{thm:rigid-W} holds
for $\mu = \delta_x$ for any $x \in X$.
Then $N \in \N$ and $( X, d ,\mathfrak{m} )$ is isomorphic
to the Euclidean space $\R^N$ with the canonical metric measure structure,
up to positive multiplicative constant on the measure.
\end{Thm}

\begin{Proof}
Since $X$ is non-compact by Lemma~\ref{lem:hat-F},
$X$ must contain more than two points.
Let $x,y \in X$.
Then Theorem~\ref{thm:rigid-W} yields that $X$ is $(0,N)$-cone
with vertex at $x$. Thus the unique minimal (unit speed) geodesic
from $x$ to $y$ can be extended to a geodesic ray.
By the same reason, the unique minimal (unit speed) geodesic
from $y$ to $x$ can be extended to a geodesic ray.
By concatenating these two geodesic rays, we obtain a line in $X$.
Thus Gigli's splitting theorem \cite{Gigli_Split} yields that
there exists an ${\sf RCD} (0,N-1)$ space $(Y, d_Y , \mathfrak{m}_Y )$
such that $( X, d, \mathfrak{m} )$ is isomorphic to $\R \times Y$ if $N \ge 2$.
When $N < 2$, $Y$ consists of one point and hence $X \simeq \R$.
In the latter case, $N = 1$ must hold
since there is $c > 0$ such that
$\mathfrak{m} ( B_r (x) ) = c r^N$ by Lemma~\ref{lem:hat-F}.
Suppose that the former case happens.
By the same reason as in the last argument,
$Y$ must contain more than two points since $N \ge 2$.
Let $x_1 , y_1 \in Y$.
Then, by applying the same argument as above to $( s, x_1 )$ and $( s, y_1 )$
instead of $x$ and $y$, we obtain a line in $X$ passing through
$( s, x_1 )$ and $( s , y_1 )$. Then we obtain the corresponding line
in $Y$ containing $x_1$ and $y_1$.
Then we can apply Gigli's splitting theorem to $Y$.
Accordingly, we can repeat the same argument to obtain the conclusion.
\end{Proof}

\begin{Rem}
We have a simpler proof of Theorem~\ref{thm:s-rigid-W} based on the notion of
tangent cones. We just give an outline of the proof here.
By \cite[Theorem~1.1]{Gigli:2013vi}, $\mathfrak{m}$-a.e. points in $X$
has a Euclidean tangent cone and the dimension $k$ of the cone satisfies $k \le N$.
We choose such a point $x \in X$.
As a consequence of Theorem~\ref{thm:rigid-W} applied to $\mu = \delta_x$,
tangent cones at $x$ of $ (X, d, \mathfrak{m})$
is unique and isomorphic to $(X, d, \mathfrak{m})$ itself
(more precisely, isomorphic to $( X, d, \mathfrak{m}, x)$
as pointed metric measure space). Thus, by the choice of $x$,
$( X, d, \mathfrak{m} )$ is isomorphic to $\R^k$ for some $k \in \N$
with $k \le N$. Then $k = N$ must hold by Lemma~\ref{lem:hat-F}.

On the one hand, we can see from this alternative proof that
it is sufficient to assume that there exists
a measurable $A \subset X$ with $\mathfrak{m} (A) > 0$ such that
the assumption of Theorem~\ref{thm:rigid-W} holds for $\mu = \delta_x$
for any $x \in A$. On the other hand, the first proof requires
essentially just \emph{$(N+1)$-points} satisfying the assumption
which are located so that ``they span the whole space''.
\end{Rem}

Finally, as a fourth result of this section,
we discuss an almost rigidity.
To state it, we borrow the notion of pointed Gromov weak distance
$\mathrm{p\mathbb{G}w}$ from \cite[Definition~3.13]{Gigli:2013wi}
between pointed metric measure spaces. In our situation below,
a convergence in $\mathrm{p\mathbb{G}w}$
(pointed measured Gromov convergence in the terminology of \cite{Gigli:2013wi})
is equivalent to a convergence in the pointed measured Gromov-Hausdorff topology
(See \cite[Theorem~3.30 and Theorem~3.33]{Gigli:2013wi}).
We choose $\mathrm{p\mathbb{G}w}$ just for simplicity of the statement.
We omit the definition of $\mathrm{p\mathbb{G}w}$ here,
but use properties of it instead.
For brevity of presentation, we state only the case $N \ge 2$,
but the corresponding assertion holds for $N \in [ 1, 2 )$.

\begin{Thm}[Almost rigidity] \label{thm:a-rigid}
Suppose $N \ge 2$.
Let $T > 0$ and $r_n : ( 0 , T ] \to ( -\infty , 0 )$ ($n \in \mathbb{N}$)
a series of non-increasing functions such that $( r_n (t) )_{n \in \mathbb{N}}$
is non-decreasing for each $t \in ( 0 , T ]$
with
\[
\lim_{t \downarrow 0} \lim_{n \to \infty} r_n ( t ) =  \sup_{n,t} r_n (t) = 0.
\]
%
Fix $s > 0$.
Let $\mathcal{M}_l$ ($l \in \mathbb{N}$) be
the set of all pointed $\mathsf{RCD} (0,N)$ spaces
$(X, d, \mathfrak{m}, x^* )$ satisfying
\begin{equation} \label{eq:W-ar0}
\mathcal{W} ( P_{t'+s} \delta_{x^*} , t'+s ) - \mathcal{W} ( P_s \delta_{x^*} , s ) \ge
r_l (t) t' \, \mbox{ for all $t' \in ( 0 , t ]$},
\end{equation}
and let $\hat{\mathcal{M}}$ be as follows:
\[
\hat{\mathcal{M}}
=
\left\{
  ( \hat{X}, \hat{d}, \hat{\mathfrak{m}} , \hat{x}^* )
  \left|
  \begin{array}{l}
      \mbox{There exists an $\mathsf{RCD}^* (N-2, N-1)$ space $( Y, d_Y , \mathfrak{m}_Y )$}
      \\
      \mbox{such that $( X, d, \mathfrak{m} )$ is $(0,N-1)$-cone built over}
      \\
      \mbox{$( Y, d_Y , \mathfrak{m}_Y )$ with vertex $\hat{x}^*$}
  \end{array}
  \right.
\right\}.
\]
%
%
Then we have
\begin{equation*}
\lim_{l \to \infty}
\sup_{
  ( X, d, \mathfrak{m}, x^*)
  \in
  \mathcal{M}_l
}
\inf_{
  ( \hat{X}, \hat{d}, \hat{\mathfrak{m}}, \hat{x}^* )
  \in
  \hat{\mathcal{M}}
}
\mathrm{p\mathbb{G}w} (
  ( X , d , \mathfrak{m}^* , x^*) ,
  ( \hat{X} , \hat{d} , \hat{\mathfrak{m}}^* , \hat{x}^* )
) = 0,
\end{equation*}
where $\mathfrak{m}^*$ (resp.\ $\hat{\mathfrak{m}}^*$)
is a normalization of $\mathfrak{m}$ (resp.\ $\hat{\mathfrak{m}}$)
so that $\mathfrak{m}^* ( B_1 (x^*) ) = 1$
(resp. $\hat{\mathfrak{m}}^* (B_1(\hat{x}^*)) = 1$).
\end{Thm}

\begin{Proof}
We first recall that the family of pointed normalized ${\sf RCD} (0,N)$ spaces
$(X,d,\mathfrak{m}^* ,x^*)$
(where ``normalized'' means $\mathfrak{m}^* (B_1 (x^*) ) = 1$)
are compact with respect to $\mathrm{p\mathbb{G}w}$.
It follows by combining
Lemma 3.32, Theorem 3.30, Remark 3.29,
comments at the beginning of Subsection 4.2
and Theorem 7.2 of \cite{Gigli:2013wi}.

Suppose that the conclusion does not hold.
Then there are $\varepsilon > 0$ and
an increasing sequence $l_n \in \N$ ($n \in \mathbb{N}$)
such that there exists $( X_n , d_n , \mathfrak{m}_n , x_n^* ) \in \mathcal{M}_{l_n}$
satisfying
\begin{align} \label{eq:W-ar1}
\inf_{
  ( \hat{X}, \hat{d}, \hat{\mathfrak{m}}, \hat{x}^* )
  \in
  \hat{\mathcal{M}}
}
\mathrm{p\mathbb{G}w} (
  ( X_n , d_n , \mathfrak{m}_n^* , x_n^*) ,
  ( \hat{X} , \hat{d} , \hat{\mathfrak{m}}^* , \hat{x}^* )
)
\ge \varepsilon
\end{align}
for each $n \in \N$.
Then there exists a convergent subsequence of
$( ( X_n , d_n , \mathfrak{m}^*_n , x^* ) )_{n \in \N}$
with respect to $\mathrm{p\mathbb{G}w}$.
We may assume that
$( ( X_n , d_n , \mathfrak{m}^*_n , x^* ) )_{n \in \N}$
itself converges without loss of generality.
We denote the limit by $( X, d, \mathfrak{m}^* , x^* )$ and the remark
at the beginning of this proof tells us that
$(X,d,\mathfrak{m}^* ,x^* )$ is a (normalized) ${\sf RCD}(0,N)$ space.
Let $P^{(n)}_t$ (resp.\ $P_t$) be the heat semigroup
on $( X_n , d_n , \mathfrak{m}_n^* )$
(resp.\ $( X, d, \mathfrak{m}^* )$).
By \cite[Theorem~7.7]{Gigli:2013wi},
$\Ent ( P_{t+s}^{(n)} \delta_{x_n^*} ) \to \Ent ( P_{t+s} \delta _{x^*} )$
for $t \in [ 0, T ]$ and
$I ( P_{t+s}^{(n)} \delta_{x^*_n} ) \to I ( P_{t+s} \delta_{x^*} )$
for a.e.\ $t \in [ 0 , T ]$. Let $J \subset [ 0 , T ]$ be the set
of points where the latter convergence occurs.
Take $t \in ( 0 , T ]$ and $t' , t'' \in J$ with $t'' < t' < t$.
Then, by Lemma~\ref{lem:mono-I}, we have
\begin{align*}
r_{l_n} (t) t'
& \le
\mathcal{W} ( P^{(n)}_{ s + t'} \delta_{x_n^*} , s + t' )
-
\mathcal{W} ( P^{(n)}_{s} \delta_{x_n^*} , s )
\\
& \le
\mathcal{W} ( P^{(n)}_{s + t'} \delta_{x_n^*} , s + t' )
- s I ( P^{(n)}_{s + t''} \delta_{x_n^*} ) + \Ent ( P^{(n)}_s \delta_{x_n^*} )
+ \frac{N}{2} \log s .
\end{align*}
Let $r_\infty := \lim_{n \to\infty} r_n$. Then, by taking $n \to \infty$
in the last inequality together with Lemma~\ref{lem:mono-I}, we have
\[
r_{\infty} (t) t'
\le
\mathcal{W} ( P_{s + t'} \delta_{x^*} , s + t' ) -
\mathcal{W} ( P_{s} \delta_{x^*} , s )
\]
for all $t' \in [ 0 , t ]$.
Since $\displaystyle \lim_{t \downarrow 0} r_\infty (t) = 0$,
the last inequality easily yields
\[
\varlimsup_{t \downarrow 0}
\frac{
  \mathcal{W} ( P_t \delta_{x_\infty^*} , t + s )
  -
  \mathcal{W} ( P_{s} \delta_{x_\infty^*} , s )
}{ t }
= 0
\]
with the aid of Theorem~\ref{thm:mono-W}.
Thus, by Theorem~\ref{thm:rigid-W},
we obtain
$
( X , d , \mathfrak{m}^* , x^* )
\in
\hat{\mathcal{M}}
$.
It contradicts with \eqref{eq:W-ar1}
and hence the conclusion follows.
\end{Proof}

\noindent
\emph{Acknowledgments.} This work was started when the authors attended the Workshop on Geometry and Probability held at the Luxembourg University in October 2013. The authors would like to thank Anton Thalmaier for his invitations which lead this work possible. The first author warmly thanks to Shouhei Honda
for his suggestion to formulate Theorem~\ref{thm:a-rigid} and
to Nicola Gigli for an improvement of Theorem~\ref{thm:a-rigid}.
He also wish to tell his gratitude to Karl-Theodor Sturm
for fruitful discussions. Especially, Theorem~\ref{thm:s-rigid-W}
comes from a discussion with him. The second author would like to thank Songzi Li for fruitful collaboration on the study of the $\mathcal{W}$-entropy for heat flows,  geodesic and Langevin flows on the Wasserstein space over manifolds. He also would like to express his gratitude to Kazuhiro Kuwae and Yu-Zhao Wang for valuable discussions on various topics related to this work.

\bibliographystyle{amsplain}
\providecommand{\bysame}{\leavevmode\hbox to3em{\hrulefill}\thinspace}
\providecommand{\MR}{\relax\ifhmode\unskip\space\fi MR }
\providecommand{\MRhref}[2]{%
  \href{http://www.ams.org/mathscinet-getitem?mr=#1}{#2}
}
\providecommand{\href}[2]{#2}

\end{document}